\documentclass[12pt]{amsart}
\usepackage{amsmath,amsthm,amscd,amsfonts,amssymb,euscript}
\swapnumbers
\usepackage[all]{xypic}
\usepackage{latexsym}
\usepackage{color} 
\theoremstyle{plain}

\begin{document}
\input{amssym.def}

\numberwithin{equation}{section}

\newtheorem{guess}{\sc Theorem}[section]
\newcommand{\bth}{\begin{guess}$\!\!\!${\bf }~}
\newcommand{\eeth}{\end{guess}}

\newtheorem{propo}[guess]{\sc Proposition}

\newcommand{\bprop}{\begin{propo}$\!\!\!${\bf }~}
\newcommand{\eprop}{\end{propo}}

\newtheorem{lema}[guess]{\sc Lemma}
\newcommand{\blem}{\begin{lema}$\!\!\!${\bf }~}
\newcommand{\elem}{\end{lema}}

\newtheorem{defe}[guess]{\sc Definition}
\newcommand{\bdefe}{\begin{defe}$\!\!\!${\it }~}
\newcommand{\edefe}{\end{defe}}

\newtheorem{coro}[guess]{\sc Corollary}
\newcommand{\bcor}{\begin{coro}$\!\!\!${\bf }~}
\newcommand{\ecor}{\end{coro}}

\newtheorem{rema}[guess]{\it Remark}
\newcommand{\brem}{\begin{rema}$\!\!\!${\it }~\rm}
\newcommand{\erem}{\end{rema}}

\theoremstyle{remark}

\newtheorem{example}[guess]{\it Example}

\newcommand{\spec}{{\rm Spec}\,}
\newtheorem{notation}{Notation}[section]
\newcommand{\bnot}{\begin{notation}$\!\!\!${\bf }~~\rm}
\newcommand{\enot}{\end{notation}}
\numberwithin{equation}{subsection} 
\newcommand{\beqa}{\begin{eqnarray}}
\newcommand{\eeqa}{\end{eqnarray}}
\newtheorem{thm}{Theorem}
\theoremstyle{definition}
\newtheorem{say}[guess]{}
\newtheorem{hint}[thm]{Hint}
\newcommand{\bsem}{\begin{say}$\!\!\!${\it }~~\rm}
\newcommand{\esem}{\end{say}}

\newsymbol \circledarrowleft 1309

\newcommand{\ha}{\mathfrak h}
\newcommand{\g}{\mathfrak g}
\newcommand{\ta}{\mathfrak t}
\newcommand{\s}{\mathfrak s}
\newcommand{\ctext}[1]{\makebox(0,0){#1}}
\setlength{\unitlength}{0.1mm}

\newcommand{\wt}{\widetilde}
\newcommand{\Lr}{\Longrightarrow}
\newcommand{\Aut}{\mbox{{\rm Aut}$\,$}}
\newcommand{\ul}{\underline}
\newcommand{\ol}{\bar}
\newcommand{\lr}{\longrightarrow}
\newcommand{\bc}{{\mathbb C}}
\newcommand{\bp}{{\mathbb P}}
\newcommand{\bz}{{\mathbb Z}}
\newcommand{\bq}{{\mathbb Q}}
\newcommand{\bn}{{\mathbb N}}
\newcommand{\bw}{{\mathbb W}}

\newcommand{\cl}{{\mathcal L}}
\newcommand{\cv}{{\mathcal V}}
\newcommand{\cf}{{\mathcal F}}
\newcommand{\cb}{{\Lambda}}
\newcommand{\mfc}{{\mathfrak C}}
\newcommand{\ce}{{\mathcal E}}
\newcommand{\co}{{\mathcal O}}
\newcommand{\cs}{{\mathcal S}}
\newcommand{\cg}{{\mathcal G}}
\newcommand{\ca}{{\mathcal A}}
\newcommand{\hra}{\hookrightarrow}
\newcommand{\mfu}{{\mathfrak U}}

\newtheorem{ack}{Acknowledgments}       
\renewcommand{\theack}{}

\title[Tensor product theorem]{Tensor product theorem for Hitchin pairs 
\\
-An algebraic approach} 
\author{V. Balaji} 
\address{Chennai Mathematical Institute SIPCOT IT Park, Siruseri-603103, India,
balaji@cmi.ac.in}
\author{A.J. Parameswaran}
\address{Kerala School of Mathematics, Kozhikode, Kerala and School of Mathematics, Tata Institute of Fundamental Research, Mumbai-400095, India, param@math.tifr.res.in}
\date{}
\thanks{The research of the first author was partially supported by the J.C. Bose Fellowship.}

\begin{abstract}
We give an algebraic approach to the study of Hitchin pairs and prove the tensor product theorem for Higgs semistable Hitchin pairs over smooth projective curves defined over algebraically closed fields $k$ of characteristic $0$   and characteristic $p$, with $p$ satisfying some natural bounds. We also prove the corresponding theorem for polystable bundles.
\end{abstract}

\maketitle
\small
\tableofcontents
\normalsize

\vspace{2mm}
\noindent

\section{Introduction} Let $X$ be a smooth projective curve over an algebraically closed field $k$. When the ground field $k$ is $\mathbb C$, the notion of a Hitchin pair is due to Nigel Hitchin. In (\cite{hitchin}, \cite{hitchin1}) he proves the basic theorem that the category of semistable Hitchin pairs of degree $0$ is equivalent to the category of complex $GL(n)$ representations of the fundamental group $\pi_1(X)$. One of the fundamental consequences of this correspondence is that the tensor product of two semistable Hitchin pairs of degree $0$ is again semistable. The Kobayashi-Hitchin correspondence in the setting of Hitchin pairs  has been generalized and extended in a number of ways starting with the far reaching one by C. Simpson (\cite{sim0},\cite{sim1}). Simpson develops the more general objects which he terms $\Lambda$--modules. This gives differential geometric proofs of the tensor product theorem for Hitchin pairs as well as for objects such as $\Lambda$--modules. In (\cite{bisschu}) Biswas and Schumacher prove simlar results for stable Higgs sheaves over arbitrary Kahler manifolds.

For the classical case of semistable bundles the tensor product theorem is usually derived as a consequence of the Narasimhan-Seshadri theorem or by using the usual Kobayashi-Hitchin correspondence. The first purely algebraic proof of the tensor product theorem is due to Bogomolov (\cite{bogo}) and little later by Gieseker(\cite{gies}) using very different methods. The third approach due to Ramanan and Ramanathan (\cite{rr}) has the advantage of being amenable to generalizations to the positive characteristic case as well (cf. \cite{imp} and \cite{bapa}). 

The aim of this paper is to give algebraic proofs of the tensor product theorem for Hitchin pairs  over ground fields of all characteristics. Towards this, we need to first develop a purely algebraic notion of Hitchin
 schemes, an object dual in a certain sense to a Hitchin pair. This is  indispensable for the algebraic proof since the standard methods of proof for usual principal bundles do not apply for the setting of principal Hitchin pairs; the Higgs structure datum has to be suitably incorporated in the algebraic setting. In the present paper, our approach, following Nori (\cite{nori}), is a Tannakian one and the notions of  ``associated" Hitchin schemes (analogous to ``associated fibrations" to principal bundles) and geometric Higgs sections arises naturally.

 We then use this new machinery for our purposes, along with  a suitable modification of the invariant theoretic ideas due to Bogomolov (\cite{bogo}) and Ramanan-Ramanathan (\cite{rr}). In positive characteristics we use the results due to Kirwan (\cite{kirwan}) and Hesselink (\cite{hesslink}) synthesized  with the methods of Ramanan-Ramanathan. This is absolutely essential in the setting of Hitchin pairs since the reduction of structure group to the Kempf-Rousseau parabolic, which is key to proof of the main theorem, is realizable geometrically only if we employ Kirwan's stratification. Representation theoretic bounds such as {\em low heights} (see \cite{imp}) come up as expected when we work in char $p$ (see Theorem \ref{maintheorem}).  In characteristic zero we generalize Bogomolov's approach to the setting of Hitchin pairs and give a different proof of the main theorem; we do this for its \ae sthetic elegance. For generalizing Bogomolov's results, we find the expos\'e due to Rousseau (\cite{rou}) just the right one and we use it freely. It would be very interesting to compare these two methods of proof since Bogomolov also provides a stratification of the unstable locus (as does Kirwan) and defines the concept of a {\em model} which is in a sense ``universal" for instability.

A word about the central principle which underlies these ``algebraic" proofs of tensor product theorems. The idea is to connect the concepts of Higgs semistability of  $G$--Hitchin pair $(E, \theta)$ (see Definition \ref{stabilitynotions}) with that of Bogomolov stability of Higgs sections of associated Hitchin pairs $(E(V), \theta_{_V})$ via a representation $G \to GL(V)$ (see Definition \ref{bogohiggs}). Since we work in the setting of Hitchin pairs, we need to work with Higgs sections of the associated objects.

The new result that emerges by this approach, apart from the \ae sthetics of a purely algebraic proof, are variations in positive characteristics for the notions of principal Hitchin pairs. We observe that the Frobenius pull-back is an inconsequential operation for Hitchin pairs and therefore notions such as {\em strong} semistability  do not provide anything new in the Higgs setting. But we show that the bounds ({\em height and separability index}) developed in \cite{imp} and \cite{bapa} are immediately applicable. In the context of the recent work of Ng\^o Bao Ch\^au  (\cite{ngo}) we believe that our approach could be of interest in positive characteristics (see also \cite{lp}).  The following theorems are the main results in the paper:

\bth (Theorem \ref{maintheorem1}, Theorem \ref{maintheorem}) Let $(V_1, \theta_1)$ and $(V_2, \theta_2)$ be two Higgs semistable Hitchin pairs with $det(V_i) \simeq \co_X, i = 1,2$. Suppose that the ground field $k$ has characteristic $p$ such that
\[
rank(V_1) + rank(V_2) < p + 2
\]
Then the tensor product $(V_1 \otimes V_2, \theta_1 \otimes 1 + 1 \otimes \theta_1)$ is also Higgs semistable. 
\eeth

\bth (Theorem \ref{polystablehiggs}) Let $(E,\theta)$ be a stable Hitchin pair of {\em degree zero} with $G$ semisimple and $\rho: G \rightarrow SL(M)$  be a  representation. Let ${\overline {\psi_{_G}(M)}}$ be as in \eqref{psibar} and Definition \ref{lowsep}. Suppose that $p > {\overline {\psi_{_G}(M)}}$. Then the associated Hitchin pair $(E(M), \theta_M)$ is polystable.
\eeth

The layout of the paper is as follows: in Section 2 we develop the generalities about Hitchin pairs and define the concept of a Hitchin scheme. In Section 3 we study principal $G$--Hitchin pairs and the associated Hitchin schemes. Section 4 contains generalization of Bogomolov stability of sections in the Hitchin pair setting. In Section 5 we recall results from the papers of Kempf, Hesselink, Kirwan and Ramanan-Ramanathan. In Section 6 and Section 7 we give an approach following Bogomolov for the proof of the main theorem in char $0$.  In Section 8 we prove the main semistability theorem in positive characteristics with the low height assumptions. In the last section we prove the theorem on polystability of associated bundles under assumptions that the characteristic $p$ is larger than the  low separability index of some natural representation spaces (see Remark \ref{correction} for some clarifications on earlier papers on this result, which treat bundles without Higgs structures). The final remarks (Remark \ref{higherdim}) indicate how these notions easily generalize to the case when $X$ is a higher dimensional variety since we work with $\mu$--semistability. 

{\it Acknowledgements}: We thank Jochen Heinloth, Madhav Nori and D.S. Nagaraj for some helpful discussions and suggestions. We sincerely thank the referee for his/her comments and suggestions. They have gone a long way to clarify the paper. The first author thanks the hospitality of TIFR and KSOM. The second author thanks the hospitality of CMI and IMSc.

\section{Hitchin pairs, basic facts}
Throughout this paper, unless otherwise stated, we have the following
notations and assumptions:

\bsem\label{uno}{\bf The Category of Hitchin pairs.}  Let $X$ be a smooth projective curve over $k$.  Let $\mathfrak{U}$ be a locally free $\co_X$--module.

\bdefe\label{higgsstructure} A {\em $\mathfrak U$--Higgs structure} (or simply a Higgs structure, since $\mathfrak U$ is fixed) on a locally free $\co_X$--module  $\cf$ is an $\co_X$--module map $\theta: \cf \to \cf \otimes \mathfrak{U}$
together with the integrability condition $\theta \wedge \theta = 0$, where  $\theta \wedge \theta:\cf \to \cf \otimes \wedge^2(\mathfrak{U})$. \edefe

\bdefe A {\em Hitchin pair}  is a locally free $\co_X$--module $\cf$ which is equipped with a Higgs structure $\theta$ and we denote it by $(\cf, \theta)$. 
\edefe

\brem The structure sheaf $\co_X$ of the base curve $X$ will always carry the {\em trivial Higgs structure}, i.e the zero map $ \co_X \to \co_X \otimes \mathfrak{U}$ unless otherwise stated.\erem

\bdefe\label{hitchinsection} The {\it space of Higgs sections} of a Hitchin pair $(\cf, \theta)$ is defined to be the space of sections $s$ of the $\co_X$--module $\cf$ such
that $\theta \circ s = 0$.\edefe

 The {\em tensor product}  of two Hitchin pairs  $ (\ce, \phi) \otimes (\cf, \theta)$ has $\ce \otimes \cf$ as the underlying bundle and the Higgs structure is defined as 
\beqa
\phi \otimes 1 + 1 \otimes \theta.
\eeqa

\bdefe\label{dualhitchin}  The {\em dual Hitchin pair is defined as the pair where $(\cf^*, -\theta^{t})$, where $\cf^*$ is the usual dual of $\cf$ and $\theta^{t}$ is defined as follows}: consider the Higgs structure $\theta: \cf \to \cf \otimes \mfu$. Taking duals, we get $\theta^*: \cf^* \otimes \mfu^* \to \cf^*$. Tensor this with $\mfu$ to get $1 \otimes \theta^* :\cf^* \otimes (\mfu \otimes \mfu^*) \to \cf^* \otimes \mfu$. Now embed $\cf^* \hra \cf^* \otimes (\mfu \otimes \mfu^*)$ using the identity section $\co_X \to {\mathcal Hom}(\mfu,\mfu) = (\mfu \otimes \mfu^*)$. Composing these maps we get
\beqa
\theta^{t}: \cf^* \to \cf^* \otimes \mfu.
\eeqa
The dual Higgs structure is given by taking $- \theta^{t}$ to be the Higgs structure on $\cf^*$. \edefe
\brem The sign $- \theta^{t}$ is given to take care that $(\cf, \theta) \otimes (\cf^*, - \theta^{t})$ gives $(\co_X, 0)$ when rank $\cf$ is 1 (cf. \cite[Page 14]{sim0}).\erem

Morphisms of Hitchin pairs are defined as usual, i.e morphisms of the $\co_X$--modules compatible with the Higgs structures. We denote by $Hitch(X)$ the category of Hitchin pairs with the tensor structure, duals and morphisms as described above. 

\brem A Higgs section can also be thought of as a Higgs morphism $s: \mathcal O_X \to \cf$, where $\mathcal O_X$ is given the trivial Higgs structure, i.e the zero map $\co_X \to \mathfrak U$.\erem

\esem

\bsem{\bf Hitchin pairs as $\cb$-modules.}
C. Simpson in (\cite[Section 2]{sim1}) gives an equivalent description of Hitchin pairs; although this is stated under assumptions of characteristic zero, it is not hard to see that the formalism holds good over positive characteristics as well. Let $\cb$ be the $\co_X$--algebra defined by 
\beqa
\cb = Sym(\mathfrak{U}^*)
\eeqa

If $\alpha \in \mathfrak {U}^*$ and if $V$ and $W$ are $\cb$--modules, then $V \otimes_{\co_X} W$ gets a $\cb$--module structure by the Leibnitz formula
\beqa\label{caution1}
\alpha(v \otimes w) = \alpha(v) \otimes w + v \otimes \alpha(w).
\eeqa
Then we have
\blem (cf. \cite[Lemma 2.13, page 85]{sim1})
\begin{enumerate}
\item Giving a Higgs structure on an $\co_X$--module $\cf$ is equivalent to giving a $\cb$--module structure on $\cf$. 

\item Morphism of Hitchin pairs are equivalently morphism of $\co_X$--modules which are simultaneously also $\cb$--module maps.

\item We have an equivalence of categories between Hitch(X) and  localy free ${\cb}$--Modules.
\end{enumerate}
\elem

\brem
In fact, Simpson (\cite[Section 2, page 77]{sim1}) considers more general objects such as bundles with integrable connections. For example, if we take ${\mathcal D}_X$ to be the sheaf of differential operators on a smooth complex curve then we could work with the category of left ${\mathcal D}_X$--modules. \erem
\esem

\bsem{\bf Hitchin algebras and Hitchin schemes.} We work with the category of affine $X$--schemes. The generalities that we develop here are essential in the paper.

\bdefe A Hitchin $\co_X$--algebra is a faithfully flat $\co_X$--algebra $\ca$ such that
\begin{enumerate}
\item $\ca$ gets a $\cb$--module structure, i.e a map
\[
\theta: \ca \to \ca \otimes_{\co_X} {\mathfrak U}
\]
\item Furthermore, for the natural $\cb$--module structure on $\ca \otimes \ca$, the multiplication map
\[
\ca \otimes_{\co_X} \ca \to \ca 
\]
and the map $\co_X \to \ca$, given by the unit in $\ca$,
are $\cb$--module maps.
\end{enumerate}

A Hitchin $X$--scheme
 is an affine $X$--scheme $f:Z \to X$ such that $f_*(\co_Z)$ gets the structure of a Hitchin $\co_X$--algebra. 
In particular, for $f = id_X$, the trivial Higgs structure on $\co_X$ gives a Hitchin $X$--scheme
 structure on $X$.

\edefe

\brem Equivalently (following Beilinson-Drinfeld (\cite{beil})), a Hitchin $\co_X$--algebra is a faithfully flat $\co_X$--algebra in the tensor category of $\cb$--modules. For example, $\co_X$ is a Hitchin $\co_X$--algebra and if $B$ is a commutative $k$--algebra, then $B \times _k \co_X$ is a Hitchin algebra. A word of caution here: a Hitchin $\co_X$--algebra is {\em not} a $\Lambda$--algebra in the usual sense of the term as can be seen from \eqref{caution1}. \erem

Let $\mfc$ denote the category of Hitchin $X$--schemes. A morphism between two   Hitchin $X$--schemes is a morphism $\phi:Z \to Y$ which preserves the $\cb$--module structure, i.e the canonical map
\beqa
\xymatrix{
Z \ar[dr]_{f} \ar[rr]^{\phi}& & Y \ar[dl]^{g} \\
& X &
}
\eeqa 
where $\phi$ induces a morphism of $\co_X$--algebras $g_*(\co_Y) \to f_*(\co_Z)$ which should be also a $\cb$--module map.

\blem\label{fibreproduct} (Fibre products in $\mfc$) Let $Z$ and $T$ be in $\mfc$. Then the fibre product $Z \times_X T$ is in $\mfc$. \elem

\noindent
\noindent {\it Proof}:
 This is clear if one uses \eqref{caution1}.
 
\bcor\label{higgsinverseimages} Let $\phi:Z \to T$ be a morphism of Hitchin $X$--schemes. Let $T_1 \subset T$ be a closed Higgs subscheme. Then the inverse image scheme $Z_1 = \phi^{-1}(T_1) \subset Z$, being a fibre product, is a closed Higgs subscheme of $Z$. \ecor

\brem Let $K = k(X)$ be the function field of $X$. We observe that we can define Hitchin algebras over $K$ as follows: let $A$ be a finite type $K$--algebra and fix a finite dimensional projective $K$--module ${\mathfrak U}_K$. Let $\Lambda_K = Sym({\mathfrak U}_K^{*})$. A Higgs structure is a map
\[
\theta: A \otimes_K {\mathfrak U}_K^* \to A
\]
\item Furthermore, for the natural $\cb_K$--module structure on $A \otimes A$, the multiplication map
\[
A \otimes_K A \to A
\]
is a $\Lambda_K$--morphism.
A Hitchin scheme
 over $K$ is $Spec(A)$ for a Hitchin algebra $A$ over $K$.
\erem

\esem

\section{Hitchin functors and principal bundles}

\bsem{\bf Hitchin functors following Nori.} Let $G$ be an affine group scheme defined over an algebraically closed field $k$. A {\em $G$--Hitchin functor}  is a tensor functor $F: Rep(G) \to Hitch(X)$ satisfying  Nori's axioms, namely $F$ is a strict, exact and faithful tensor functor (cf. \cite[Page 77]{nori}) such that the following diagram commutes:
\beqa
\xymatrix{
                  &Hitch(X) \ar[d]^{``forget"} \\
Rep(G) \ar[ru]^{F} \ar[r]^{F'} & Vect(X)  \\
}
\eeqa
where the functor $forget: Hitch(X) \rightarrow Vect(X)$  forgets the Higgs structure. If $V$ is a finite dimensional $G$--module, we will denote the associated Hitchin pair by $F(V)$. Note that the data underlying $F(V)$ is a locally free sheaf $F'(V)$ together with a Higgs structure on $F'(V)$.

\brem Observe that the {\em forget} functor is a tensor functor in this situation and hence by Nori's observation, the functor $F'$ canonically gives rise to a principal $G$--bundle on $X$. The aim in this section is to represent the functor $F$ by a suitable ``Hitchin scheme" which has as its underlying principal $G$--bundle the one given by $F'$. \erem

Let ${\mathcal S}(X)$ be the category of quasi-coherent $\co_X$--modules. We extend $F$ to a functor $$\bar {F} :\{G-mod\} \to {\mathcal S}(X)$$ as follows:
Let $M$ be an arbitrary $G$--module. Express $M$ as a direct limit of finite dimensional $G$--modules 
\beqa
M = \underset{j}{\underset{\longrightarrow}{\lim}}~V_j
\eeqa
now define,  
\beqa
{\bar{F}}(M) := \underset{j}{\underset{\longrightarrow}{\lim}}~F(V_j)
\eeqa
which realises $ {\bar{F}}(M)$ as a quasi-coherent $\co_X$--module. 
\brem The definition of the quasi-coherent sheaf $ {\bar{F}}(M)$ is independent of the particular limit chosen. \erem
\esem

\bsem{\bf Conjugate Higgs structure.} We observe that $ {\bar{F}}(M)$ has a natural Higgs structure as follows: 
for each $V_j$ we have the Higgs structure given by 
\beqa\label{psi}
\psi_j:F(V_j)  \to F(V_j) \otimes \mathfrak{U}
\eeqa
Taking limits we get the Higgs structure on $ {\bar{F}}(M)$. In this situation, we equip $\bar{F}(M)$ with the {\em conjugate Higgs structure} as follows:

The Higgs structure \eqref{psi} canonically induces on the dual locally free $\co_X$--module $F(V_j^*)$ and a Higgs structure $\hat{\psi_j}: F(V_j^*) \to F(V_j^*) \otimes \mathfrak{U}$. 
Dualizing $\hat{\psi_j}$ we get
\beqa\label{dualofdual}
(\hat{\psi_j})^{*} := \varphi_j:F(V_j) \otimes \mathfrak{U}^* \to F(V_j)
\eeqa

Now taking limits and observing that tensor products commutes with direct limits, we get a map
\beqa
\underset{j}{\underset{\longrightarrow}{\lim}}~{\varphi_j} = \varphi: {\bar{F}}(M) \otimes \mathfrak{U}^* \to {\bar{F}}(M)
\eeqa
which we term the {\em conjugate Higgs structure} on the quasi-coherent module ${\bar{F}}(M)$. We observe that the induced $\cb$--module structure on ${\bar{F}}(M)$ comes from this conjugate Higgs structure and extending it to an action of $Sym({\mathfrak U}^*)$.\esem

\brem We note that in the finite dimensional setting, a ``conjugate Higgs structure" is in reality the ``dual" of the ``dual Higgs structure" as defined in Definition \ref{dualhitchin}. Note the importance of the signs. This can be seen above in \eqref{dualofdual} above, where we take the ``dual Higgs structure"  $\hat{\psi_j}$ on $F(V_j^*)$ and then once more dualize to get $\varphi_j$ which is the ``conjugate Higgs structure" on $F(V_j)$. Taking ``duals" works fine in the finite dimensional setting but since we need the infinite dimensonal setting, we need to be careful. Finally, the $\Lambda$--module structure on a Hitchin pair comes via the conjugate Higgs structure in this sense. \erem

\bsem{\bf Associated Hitchin scheme.}
We have a natural extension of the functor $F$ to the category of affine $G$--schemes:
$${\mathcal H}_F: \{affine~ G-schemes\} \to \{affine ~X-schemes\}$$
To see this, let $Z = Spec(k[Z])$. Then $k[Z]$ is a (possibly) infinite dimensional $G$--module and also the multiplication map $k[Z] \otimes k[Z] \to k[Z]$ along with the tensor axiom for $F$ give the sheaf ${\bar{F}}(k[Z])$ the structure of an $\co_X$--algebra. This defines a $X$--scheme $j:{\mathcal H}_F(Z) \to X$ where 
\beqa
{\mathcal H}_F(Z):= {\mathcal{S}pec}({\bar{F}}(k[Z]) 
\eeqa
and further we can identify $j_*(\co_{{\mathcal H}_F(Z)}) = {\bar{F}}(k[Z])$ as an $\co_X$--algebra.\esem

Let $F$ be a Hitchin functor and let $Z = Spec(k[Z])$ be an affine $G$--scheme. Then since $j_*(\co_{{\mathcal H}_F(Z)}) = {\bar{F}}(k[Z])$ we have the canonical {\em conjugate Higgs structure} on the $\co_X$--algebra ${\bar{F}}(k[Z])$, (i.e a $\co_X$--module morphism)
\beqa
\eta_Z: {\bar{F}}(k[Z]) \otimes \mathfrak{U}^* \to {\bar{F}}(k[Z])
\eeqa
Again, the multiplication map $k[Z] \otimes k[Z] \to k[Z]$ along with the tensor axiom for $F$ give the $\co_X$--algebra ${\bar{F}}(k[Z])$ the structure of an Hitchin $\co_X$--algebra. 
This therefore gives the structure of a  {\em Hitchin $X$--scheme
} on ${\mathcal H}_F(Z)$,
which we denote by $({\mathcal H}_F(Z), \eta_Z)$.

\bdefe The Hitchin $X$--scheme
 $({\mathcal H}_F(Z), \eta_Z)$ is called the associated Hitchin
 scheme
 to the Hitchin functor $F$. \edefe

\bsem {\bf Geometric Higgs section.} Recall that a section of the fibration ${\mathcal H}_F(Z)$  is a $X$--morphism $s: X \to {\mathcal H}_F(Z)$ which is given by an {\em $\co_X$--algebra morphism} $ s: {\bar{F}}(k[Z]) \to \co_X$. 

\bdefe We say that $s:X \to ({\mathcal H}_F(Z), \eta_Z)$ is a {\em geometric Higgs section} of the associated Hitchin
 scheme
 $({\mathcal H}_F(Z), \eta)$ if further the composite:
\beqa
s \circ \eta: {\bar{F}}(k[Z]) \otimes \mathfrak{U}^* \to {\bar{F}}(k[Z]) \to \co_X
\eeqa
is zero. 
\edefe

\brem Equivalently, $ s: {\bar{F}}(k[Z]) \to \co_X$ is a Hitchin $\co_X$--algebra morphism. In other words, $s$ is a section in the category $\mfc$. \erem

We will denote by $\ce(V)$ the image $F_E(V)$ as a locally free $\co_X$--module and $(\ce(V), \theta_V)$ the associated Hitchin pair. While viewing the $G$-module $V$ as an affine scheme we will use the notation 
\beqa\label{mathbbv}
{\mathbb V} = Spec(Sym(V^*))
\eeqa 
The associated geometric fibre space is denoted by $E({\mathbb V})$.\esem

\brem Recall that $E(\mathbb V)$ is the {\em geometric vector bundle} in the sense of Grothendieck.\erem

\bprop\label{usualversusgeometric} A Higgs section (see Definition \ref{hitchinsection})  of the Hitchin pair $(\ce(V), \theta_V)$ gives a geometric Higgs section of the associated Hitchin
 scheme
 $(E(\mathbb V), \eta_{\mathbb V})$ and conversely. 
\eprop

\noindent {\it Proof}:
By the functorial property of the symmetric algebra, an $\co_X$--module map $\co_X \to \ce(V)$ canonically gives rise to an $\co_X$--algebra map $Sym(\ce(V^*)) \to \co_X$ i.e an $\co_X$--algebra 
map ${\bar{F}}(k[{\mathbb V}]) \to \co_X$ and conversely.

We need only observe that the Higgs section property is also preserved. But this can be formulated as $\cb$--structures and morphisms which preserve this structure. Therefore, by the functorial property of the symmetric algebra, a $\cb$--module map $\co_X \to \ce(V)$ canonically gives rise to a Hitchin $\co_X$--algebra map $Sym(\ce(V^*)) \to \co_X$ i.e a Hitchin $\co_X$--algebra
map ${\bar{F}}(k[\mathbb V]) \to \co_X$ and conversely. This takes care of the Higgs property. \begin{flushright} {\it QED} \end{flushright}

\bsem{\bf Principal Hitchin pairs.} We now define principal Hitchin pairs and show the representability of a $G$--Hitchin functor by a principal Hitchin pair. 
\bdefe A {\em  principal $G$--Hitchin pair} $j:E \to X$ is a principal $G$--bundle together with the structure of an associated Hitchin
 scheme
  on $E$, i.e a conjugate Higgs structure on the $\co_X$--algebra $j_*(\co_{E})$, viz 
\beqa\label{hitchinpaireta}
\eta: j_*(\co_{E}) \otimes \mathfrak{U}^* \to j_*(\co_{E})
\eeqa
Furthermore, $j_*(\co_{E})$ gets the structure of a Hitchin $\co_X$--algebra. In other words, $E \to X$ is a principal $G$--object in the category $\mfc$.\edefe

Denote the principal $G$--Hitchin pair by the pair $(E,\eta)$.

\brem\label{upperstr} Giving a principal $G$--Hitchin pair $(E,\eta)$ gives the structure sheaf $\co_E$ of the underlying scheme $E$ a structure of $j^*(\Lambda)$--module. This comes by taking $j^*(\eta)$ for the $\eta$ in \eqref{hitchinpaireta}. This firstly gives a Higgs structure for the locally free sheaf $j^*(\mathfrak{U}^*)$ which then extends to a $j^*(\Lambda)$ structure on $\co_E$. \erem

\bth\label{norihiggs} A principal Hitchin pair $E$ canonically defines a Hitchin functor
\beqa
F_E: Rep(G) \to Hitch(X)
\eeqa

Conversely, let $F$ be a Hitchin functor. Then there exists a principal Hitchin pair $E$ unique upto unique isomorphism such that there is a tensor equivalence of functors $F \simeq F_E$. \eeth

\noindent
{\it Proof}: In \cite[Proposition 2.9]{nori}, Nori proves this theorem for a tensor functor $F: Rep(G) \to Vect(X)$ i.e without the Higgs structures. Observe that we have an equivalence of categories $Hitch(X) \simeq \{\Lambda-mod\}$. The point to note is that taking direct limits commutes with the $\Lambda$-module structure. The representing torsor $E$ in the category of principal bundles (obtained in Nori's theorem) has the property that the associated locally free $\co_X$--modules $\ce(V)$ are $\Lambda$--modules as well. Further, the tensor structure on $Image(F_E)$ coupled with the Higgs structure gives ${\overline F}_E(k[G])$ or equivalently $j_*(\co_E)$, the structure of a Hitchin $\co_X$--algebra and we are done. 

Conversely, let $E$ be a principal $G$--Hitchin pair. Let $V \in Rep(G)$ be a finite dimensional $G$--module. Then we need to show that the associated vector bundle ${\mathcal E}(V)$ is a Hitchin pair and this association is functorial.

By \cite[Lemma 2.6]{nori}, we have a functorial isomorphism of $G$--sheaves:
\beqa
j^*({\mathcal E}(V)) \simeq V_E = V \otimes_k \co_E
\eeqa

The trivial sheaf $V_E$ gets the obvious structure of a $j^*(\cb)$--module on $\co_E$ (see Remark \ref{upperstr}). Let $\tilde{\theta}: j^*({\mathcal E}(V)) \otimes j^*(\mathfrak U^*)\to j^*({\mathcal E}(V)) $ be the induced conjugate Higgs structure. Then by the projection formula, this structure descends the structure of a $\cb$--module on ${\mathcal E}(V)$. This $\cb$--module structure is clearly functorial and proves the converse.
\begin{flushright} {\it QED} \end{flushright}

\brem Let $G$ be a connected semisimple algebraic group. We recall that when the ground field is the field of complex numbers, C. Simpson has defined a principal Hitchin pair, or a principal Higgs bundle as a principal $G$--bundle together with a section $\theta \in H^0(E(\mathfrak g) \otimes \Omega_X^1)$ with the integrability conditions. We remark that this definition is equivalent to giving a $G$--Hitchin functor and this can be seen as follows. 

Let $\mathfrak g$ be a semisimple Lie algebra over $\bc$. Then the category of Lie algebra modules $\mathfrak g \to {\mathfrak gl}(V)$ is a neutral Tannaka category and can be seen to recover back the group $G$. This is in a sense the ``infinitesimal Tannakian construction" as done for example in \cite[Proposition 6.11]{milne}. 

We carry over this formalism to the setting of Hitchin functors. Given a Hitchin functor $F: Rep(G) \to Hitch(X)$ let $\rho: G \to GL(V)$ be an object in $Rep(G)$. Then for every $V$, we have a Hitchin pair
\[
\theta_V: \co_X \to F(V)^* \otimes F(V) \otimes \mathfrak U
\]
or equivalently a section $\theta_V \in H^0(F({\mathfrak gl}(V)) \otimes \mathfrak U)$. These sections have the naturality with respect to the tensor structure on $Rep(G)$ and by the ``infinitesimal picture" mentioned above, we get the required section $\theta \in H^0(F(\mathfrak g) \otimes \mathfrak U))$.

Conversely, given $(E, \theta)$ as in Simpson, for every $\rho: G \to GL(V)$ consider the induced differential $d \rho: \mathfrak g \to {\mathfrak gl}(V)$. This induces an $\co_X$--module map:
\[
\theta: \co_X \to E(\mathfrak g) \otimes \mathfrak U \to (E({\mathfrak gl}(V)) \otimes \mathfrak U))
\]
which gives $\theta_V \in H^0(E({\mathfrak gl}(V)) \otimes \mathfrak U))$ or equivalently a Hitchin functor.\erem

\brem From now on because of the equivalence $F \simeq F_E$, we will denote the {\em associated Hitchin
 scheme
} by $(E(Z), \eta_Z)$ and the sheaf ${\bar{F}}(k[Z])$ simply by $E(k[Z])$.\erem

\bsem{\bf Associated maps.} Let $(E, \eta)$ be a $G$--Hitchin pair. A geometric Higgs section $s:X \to (E(Z), \eta_Z)$ of the associated Hitchin
 scheme
 can therefore be  viewed as a map $$t:E(k[Z])) \to \co_X$$ of $\co_X$--algebras   such that the composite
\beqa\label{higgssection}
t \circ \eta_Z: E(k[Z])) \otimes \mathfrak{U}^* \to E(k[Z])) \to \co_X 
\eeqa
is zero. {\em Note that we give as always $\co_X$ the trivial Higgs structure}.

\brem On the algebra  $k[Z]$ the conjugate Higgs structure is nothing but the Higgs structure on the {\em restricted dual} of $k[Z]$, viz, taking the conjugate Higgs structure on the finite dimensional modules and taking the limit of the duals gives the restricted dual. \erem 
\esem

We have the following central fact:
\bprop\label{assocmorphs} Let $(E, \eta)$ be a $G$--Hitchin pair. Let $Z$ and $T$ be two affine $G$--schemes and let $\phi:Z \to T$ be a $G$--map. Then $\phi$--induces a map of associated Hitchin schemes $E(\phi): (E(Z), \eta_Z) \to (E(T), \eta_T)$. Further, a geometric Higgs section $s:X \to  (E(Z), \eta)$ gets mapped to a geometric Higgs section $E(\phi) \circ s : X \to  (E(T), \eta_T)$. 
\eprop
\noindent
\noindent {\it Proof}:
 The map $\phi$ induces a map of $G$--modules
\beqa
\phi^*: k[T] \to k[Z]
\eeqa
Now express $k[T]$ as $k[T] = \underset{j}{\underset{\longrightarrow}{\lim}}~V_j$, where $V_j$ are finite dimensional $G$--modules. Similarly, $k[Z] = \underset{l}{\underset{\longrightarrow}{\lim}}~W_l $. A $G$--module map $\phi^*$ is therefore  the data which gives for every $j$ a $W_{\phi(j)}$ together with a family of $G$--module maps of finite dimensional modules
\beqa
\phi_j:V_j \to W_{\phi(j)}
\eeqa
inducing maps of the dual structures on the bundles and associated morphisms of the dual Hitchin pairs 
\beqa
\begin{CD}
  \xymatrix{ E(V_j) \otimes \mathfrak{U}^* \ar[r]\ar[d]  & E(V_j) \ar[d] 
    \\
     E(W_{\phi(j)}) \otimes \mathfrak{U}^* \ar[r] & E(W_{\phi(j)})  &\\
} 
\end{CD}
\eeqa 
Now taking limits we get
\beqa\label{assocmorph}
\begin{CD}
  \xymatrix{ E(k[T]) \otimes \mathfrak{U}^* \ar[r]\ar[d]  & E(k[T]) \ar[d] 
    \\
     E(k[Z]) \otimes \mathfrak{U}^* \ar[r] & E(k[Z])  &\\
} 
\end{CD}
\eeqa 

The vertical arrow morphism of $\co_X$--algebras
\beqa
E(k[T]) \to E(k[Z])
\eeqa
induces at the scheme level the morphism
\beqa
E(\phi):E(Z) \to E(T)
\eeqa
and the diagram of sheaves above gives that this is a morphism of associated Hitchin schemes. This proves the first part of the proposition.

By \eqref{higgssection} a Higgs section of $E(k[Z])$ is a map $E(k[Z]) \to \co_X$ such that the composite
\beqa
E(k[Z]) \otimes \mathfrak{U}^* \to E(k[Z]) \to \co_X
\eeqa
is zero. Hence by \eqref{assocmorph} we get a diagram

\beqa
\begin{CD}
  \xymatrix{ E(k[T]) \otimes \mathfrak{U}^* \ar[r] \ar[d]  & E(k[T])  \ar[d]\ar[r] & \co_X
    \\
     E(k[Z]) \otimes \mathfrak{U}^* \ar[r] & E(k[Z]) \ar[r] & \co_X \ar@{=}[u] \\
} 
\end{CD}
\eeqa 
i.e, the section $E(k[T]) \to \co_X$ is induced by the composite $E(k[T]) \to E(k[Z]) \to \co_X$.

The commutation of the left half of the diagram immediately implies that the composite
\beqa
E(k[T]) \otimes \mathfrak{U}^* \to E(k[T]) \to \co_X
\eeqa
is also zero which implies that the $E(\phi)(s)$ is also a geometric Higgs section of the associated Hitchin
 scheme
 $E(T)$.
\begin{flushright} {\it QED} \end{flushright}

\esem

\section{Bogomolov stability of sections}
Let $k$ be an algebraically closed field of arbitrary characteristic and $G$ a {\em connected reductive algebraic group}. 
\bsem\label{reductionofstrgrp}{\bf Higgs reduction of structure group.} Let $H \subset G$ be a closed subgroup. 

\bdefe\label{higgsreductionfunctorial} A reduction of structure group of a $G$--Hitchin functor is a factoring of $F$ as follows:

\beqa
\xymatrix{
Rep(G) \ar[d]_i \ar[r]^{F} & Hitch(X) \\
Rep(H)  \ar[ru]_{F'} 
}
\eeqa
where $F': Rep(H) \to Hitch(X)$ is a $H$--Hitchin functor.\edefe

\blem\label{higgsred1} Let $H \subset G$ be a closed subgroup and let $(E,\theta)$ be a principal Hitchin pair. Let $E_{_H} \subset E$ be a reduction of structure group of the underlying $G$--bundle $E$ to $H$. Suppose that $E_{_H}$ gets the structure of a Hitchin $X$--scheme
 and the inclusion $E_{_H} \hookrightarrow E$ is a morphism of Hitchin schemes. Then $E_{_H}$ is a Higgs reduction of structure group. \elem

\noindent
\noindent {\it Proof}:
 The proof is formal and follows easily from Theorem \ref{norihiggs}.

\blem\label{onemorehiggsred} Let $(E,\eta)$ be a principal Hitchin pair giving rise to  $F_E:Rep(G) \to Hitch(X)$. Let $H \subset G$ be a closed reductive subgroup. Giving a Higgs reduction of structure group of $F_E$ to $H$ is equivalent to giving a Higgs section of the associated Hitchin
 scheme
 $E(G/H) \to X$. \elem

\noindent
\noindent {\it Proof}: Let $Z$ be an affine $G$--scheme.
  A geometric Higgs section $s:X \to (E(Z), \eta_Z)$ of the associated Hitchin scheme can therefore be  viewed as a $G$--equivariant Higgs $X$--morphism: 
\beqa
\xymatrix{
E \ar[dr]_{} \ar[rr]^{\phi}& & Z \times X \ar[dl]^{pr_2} \\
& X &
}
\eeqa

\noindent
of Hitchin $X$--schemes $E$ and $Z \times X$ where the  Hitchin scheme structure on the product is the one induced from the canonical structure on $X$. Recall that $X$ always carries the ``trivial Higgs structure".

From this, it follows that any subset $T \subset Z \times X$ flat over $X$ is a Higgs subscheme. Hence, for any $G$--equivariant map from the Hitchin $G$--scheme $E$ to $Z \times X$, the inverse image will get the structure of a Higgs $X$--subscheme on $E$ which is $G$--invariant. This follows from Corollary \ref{higgsinverseimages}.

Therefore if $z \in Z$ is any point then we can consider the closed Higgs subscheme $\{z\} \times X \hookrightarrow Z \times X$. This is a Higgs subscheme for any point $z \in Z$ since the Higgs structure on $Z \times X$ is the one induced by the structure on $X$. Then by Corollary \ref{higgsinverseimages}, the inverse image subscheme $\phi^{-1}(\{z\} \times X) = E_z$ is a Higgs $X$--subscheme of $E$. 

Specializing to the case when $Z = G/H$ which is assumed affine, we  then see that a geometric Higgs section $s:X \to E(G/H)$ is given by a $G$--equivariant Higgs $X$--morphism $\phi: E \to G/H \times X$. The inverse image of the identity coset $\phi^{-1}(e.H \times X)$ then gives $E_{_H} \subset E$ as a Higgs subscheme. By classical geometry (eg. Kobayashi $\&$ Nomizu), one knows that $E_{_H} \subset E$ gives the $H$--reduction associated to $s$, hence by Lemma \ref{higgsred1} the induced $H$--reduction $E_{_H}$ is in fact a Higgs reduction of structure group of the principal Hitchin pair $(E, \eta)$ to $H$. The converse is easy to see. Again classical geometry shows that giving $E_H$ gives rise to a $G$--equivariant $X$--morphism $E \to G/H \times X$. $E$ gets a Higgs structure and this map is trivially a Higgs morphism as seen above.\begin{flushright} {\it QED} \end{flushright}

\brem\label{generichiggsred} More generally, suppose that $H \subset G$ is a subgroup such that $G/H \subset Z$ is an arbitrary subscheme of the affine $G$--scheme $Z$. Let $K = k(X)$ be the function field of the base curve $X$. Suppose further that over the generic point $\xi \in X(K)$, the reduction section $s(\xi)$ lies in $ E(G/H) \subset E(Z)$. 

Now consider the induced morphism 
\[
\phi_K: E_K \to G/H \times Spec(K) \subset Z \times Spec(K)
\]
we see that $\phi_K^{-1}(e.H \times Spec(K))$ gives a Higgs $K$--subscheme $(E_{_H})_K \subset E_K$, i.e a generic Higgs reduction of structure group to $H$. \erem

\brem\label{caution} We fuss here about the affineness of $Z$ since we have developed the earlier formalism of associated Hitchin spaces only for affine $G$--schemes. Possibly, a graded version of this would allow us projective $G$--schemes as well. In any case, when we need to talk of Higgs reductions to parabolic subgroups, we  exercise caution while interpreting the reduction datum as sections. \erem

We now make a few remarks on the compatibility of the Higgs structure with the reduction of structure group.

\blem\label{genhiggs} Let $H \subset G$ and let $E_{_H} \subset E$ be a reduction of structure group to $H$. Suppose that for a dense open $U \subset X$, the Higgs structure on $E$ comes from $E_{_H}$. Then the Higgs structure on $E$ on the whole of $X$ comes from $E_{_H}$. \elem
\noindent
\noindent {\it Proof}:
 By \cite[Proposition 2.21]{dm}, any $H$--module $W$ as an $H$--module is a subquotient of a $G$--module $V$. i.e, there is a finite dimensional $H$--submodule $M \hra V$ and a $H$--module surjection $M \twoheadrightarrow W$. So we have a diagram of $H$--modules:
\beqa
\xymatrix{
 M \ar[d] \ar[r] & V  \\
 W   
}
\eeqa

Applying the functor $E_{_H}$, we get the diagram of vector bundles:

\beqa
\xymatrix{
 E_{_H}(M) \ar[d] \ar[r] & E_{_H}(V)  \\
 E_{_H}(W)   
}
\eeqa
where $E_{_H}(V)$ has a Higgs structure since $E_{_H}(V) \simeq E_G(V)$ and over a dense open $U \subset X$, the diagram is one of Hitchin pairs. Since $E_{_H}(M) \subset E_{_H}(V)$, it follows that the $X$--Higgs structure on $E_{_H}(V)$ restricts to a $X$--Higgs structure on $E_{_H}(M)$ extending the given one on $U$. Similarly, by considering the surjection $E_{_H}(M) \to E_{_H}(W)$ we get an extension of the $U$--Higgs structure on $E_{_H}(W)$ to the whole of $X$. This implies that $E_{_H}$ gives a $H$--Hitchin functor $F_H: Rep(H) \to Hitch(X)$. \begin{flushright} {\it QED} \end{flushright}

\esem

\bsem{\bf Semistability of principal Hitchin pairs.}
We now give the definitions of Higgs semistable (resp. polystable, stable) principal Hitchin pairs. 
\esem

\bdefe\label{parabolicreductions} Let $(E, \theta)$ be a principal $G$--Hitchin pair.  A reduction of structure group $\sigma:X \to E(G/P)$ of the underlying principal $G$--bundle to a parabolic subgroup $P \subset G$ is said to be a Higgs reduction if the $P$--subbundle $E_{_P} \subset E$ (induced by $\sigma$) gives a Higgs reduction of structure group in the sense of Definition \ref{higgsreductionfunctorial}. In other words, there is a Higgs structure $\theta_{_P}$ on $E_{_P}$, such that the extension of structure groups takes $(E_{_P},\theta_{_P})$ to $(E, \theta)$.  \edefe

\brem\label{affinebusiness} Note that the above definition allows us to handle reduction of structure groups to parabolic subgroups as well. More precisely, the reduction section $\sigma$ above gives $E_{_P}$ and we impose the condition that this $P$--bundle gets the structure of a Hitchin scheme and the induced Hitchin scheme structure on $E$ is the one coming from the original Hitchin pair structure $(E, \theta)$ (see Remark \ref{caution} and the discussions before the remark). \erem

\begin{example}\label{eg} Let $(E,\theta)$ be a $G$--Hitchin pair and let $G/P \simeq {\mathbb P}(V)$ for a finite dimensional $G$--module $V$. Then giving a Higgs reduction $s:X \to E(G/P) \simeq E({\mathbb P}(V))$ is equivalent to giving a Higgs line subbundle $L \subset E(V)$ for the locally free Hitchin pair $(E(V), \theta_V)$ in the sense that, there is a Higgs structure $L \to L \otimes {\mathfrak U}$ such the inclusion $L \hra E(V)$ preserves the Higgs structures. \end{example}

\bdefe\label{stabilitynotions} (Following A. Ramanathan) 
We follow the convention that if $\chi$ is a dominant character on a parabolic subgroup $P \subset G$, then the {\em dual} $L_{\chi}^{\vee}$,  of associated line bundle $L_{\chi}$ is {\em ample}. 
\begin{enumerate}
\item\label{higgssemistable} The $G$--Hitchin pair $(E, \theta)$ is called {\em Higgs semistable} (resp. {\em Higgs stable}) if for
every parabolic subgroup $P$ of $G$, and for every Higgs-reduction of structure
group $\sigma_{_P} : X \to E(G/P)$ to $P$ and for any dominant character
$\chi$ of $P$, the bundle $\sigma_{_P}^* (L_\chi))$ has degree $\leq 0$
(resp.$<0$). Observe that, if $E_{_P}$ is the induced $P$--Hitchin scheme coming from $\sigma_{_P}$, then we have an isomorphism of  line bundles $E_{_P}(\chi) \simeq \sigma_{_P}^* (L_\chi))$ on $X$.

\item \label{admis} A Higgs-reduction of structure group of $(E, \theta)$ to a
parabolic subgroup $P$ is called {\em admissible} if for any character
$\chi$ on $P$ which is trivial on the center of $G$, the line bundle
$E_{_P}(\chi)$ associated to the reduced $P$-bundle $E_{_P}$ has degree zero.

\item\label{poly} A $G$--Hitchin pair $(E,\theta)$ is said to be {\em Higgs polystable} if it is semistable and furthermore, for every admissible reduction of structure group $(E_{_P}, \theta_{_P})$ to a parabolic subgroup $P$, there is a Levi subgroup $R \subset P$ together with a 
Higgs-reduction of structure group $(E_{_R}, \theta_{_R})$ to $R$.

\end{enumerate}
\edefe

\brem Recall the usual notions of Higgs semistability of locally free Hitchin pairs (cf. Simpson \cite{sim0}). This is analogous to the $\mu$--semistability definition, namely, $(W, \eta)$ is Higgs semistable if for every Hitchin subpair $(W_1, \eta_1)$, we have $\mu(W_1) \leq \mu(W)$, where $\mu(W) = {{deg(W)} \over {rank(W)}}$. It is the usual exercise  to show that if $(W, \eta)$ is a Higgs semistable (resp. Higgs stable, Higgs polystable) locally free Hitchin pair, then the underlying principal $GL(r)$--Hitchin pair is Higgs semistable (resp. Higgs stable, Higgs polystable) in the above sense. \erem

\brem\label{noadmis} If a $G$--Hitchin pair has no admissible reduction to a proper parabolic subgroup $P$ then the Hitchin pair is easily seen to be {\em Higgs stable}. \erem


\bdefe\label{bogohiggs}(Following Bogomolov) Let $(E,\theta)$ be a principal $G$--Hitchin pair
and let $G \lr GL(V)$ be a representation of $G$. Let $s$
be a Higgs section of the associated Hitchin
 scheme
 $(E({\mathbb V}), \theta_{\mathbb V})$. Then we call the  section
$``s"$ {\em Bogomolov stable (resp. Bogomolov semistable, Bogomolov unstable)} relative to $G$  if at one point
$x \in X$ the value of the section
$s(x)$ is {\em stable (resp semistable, unstable)} in the GIT sense, i.e as points on $\mathbb V$.\edefe 

\brem Recall from GIT the definitions of semistability, stability and instability of points of $\mathbb V$. A point $\xi \in \mathbb V$ is semistable  if $0 \notin \overline{O(\xi)}$; the point $\xi$ is stable if furthermore, $O(\xi)$ is  closed and $Stab(\xi)$ is finite. The point $\xi$ is unstable if it is not semistable.  \erem
\noindent
\brem\label{indepofpoint} It is easy to see the non-dependence of the definition on the point $x
\in X$. Consider the inclusion $k[V]^{G} \hookrightarrow k[V]$ and the
induced morphism $q:{\mathbb V} \to {\mathbb V}/G$. This induces a morphism $E(q):E({\mathbb V}) \to E({\mathbb V}/G)$. Observe that ${\mathbb V}/G = Spec(k[V]^{G})$ is a trivial $G$-space. Thus we have the
following diagram:
\beqa
s:X \lr E({\mathbb V}) \to E({\mathbb V}/G) \simeq X \times {\mathbb V}/G
\eeqa
Composing with the second projection we get a morphism $X \to {\mathbb V}/G$
which is constant by the {\em projectivity} of $X$. Hence the value of the
section is determined by one point in its $G$-orbit.  (cf. \cite[1.10]{rou}). Thus the fibre of $q$ containing the orbit $Orb(s(x))$ is independent of $x \in X$.

By GIT, one knows that an orbit $O$ consists of points which are {\em unstable}  if and only if $O \subset q^{-1}(q(0))$. Similarly, $O$ consists of {\em stable} points if and only if $O = q^{-1}(q(O))$ and furthermore, the stabilizer of a point of $O$ is finite.

Thus, the property of whether a section is Bogomolov semistable, stable or unstable is reduced to checking it at any point of the base space $X$. \erem

\blem\label{degandss} Let $G$ be a reductive group and let $E$ be a principal $G$--bundle. Let $V$ be a finite dimensional $G$--module and let $s:X \to E(V)$ be a Bogomolov semistable section of $E(V)$, or equivalently, $s(x)$ is GIT semistable for some $x \in X$. Then $deg(s^*(L)) \geq 0$ for the $G$--linearized {\em ample bundle} $L$ on ${\mathbb P}(V)$. More generally, let $Y$ be a projective variety on which there is a $G$--linearized action with respect to an ample line bundle $L$. If $s:X \to E(Y)$ is a section such that for some point $x \in X$ image $s(x) \in E(Y)_x$ is GIT semistable for the $G$--action. Then, $deg(s^*(L)) \geq 0$.
If moreover, $deg(s^*(L)) = 0$, then the section takes it values in the GIT fibre $F \subset Y^{ss} \to Y/G$ containing $s(x)$.

\elem

{\it Proof}: Since $s(x)$ is GIT semistable for the $G$--action on $V$, there exists a $G$--homogeneous polynomial of degree $n > 0$, which is non-zero on $s(x)$. In other words, we get a non-zero section of $s^*((L)^n)$, implying that $deg(s^*(L)) \geq 0$. The second half is similar and for details see \cite[Proposition 3.10]{rr}. \begin{flushright} {\it QED} \end{flushright}

\section{A summary of results on instability}

\bsem{\bf Some notations and preliminary definitions.} Let $k$ be an algebraically closed field of arbitrary characteristic. Let
$G$ be a connected reductive algebraic group over $k$.  Let $T$ be a
maximal torus of $G$ which we fix throughout and we fix a Borel subgroup $B \supset T$. Let $X(T):=Hom(T,{\bf G}_m)$ be the character group of $T$
and $Y(T):=Hom({\bf G}_m,T)$ be the 1-parameter subgroups of $T$ (defined over $k$). Let
$R\subset X(T)$ be the root system of $G$ with respect to $T$. Let
${\mathcal W}$ be the Weyl group of the root system $R$.  Let $(~,~)$ denote
the ${\mathcal W}$-invariant inner product on $X(T)\otimes {\bf R}$. This inner product determines an inner product for any other maximal torus since any two are conjugate.\esem

For 
$\alpha\in R$, the corresponding co-root $\alpha^{\vee}$ is
$2\alpha/(\alpha,\alpha)$.  Let $R^{\vee}\subset{X(T)\otimes {\bf
R}}$ be the set of all co-roots. Let $B\subset G$ be a Borel subgroup
containing $T$. This choice defines a base $\Delta^+$ of $R$ called
the {\it simple roots}. Let $\Delta^-=-\Delta^+$. A root in $R$ is
said to be {\em positive} if it is a non-negative linear combination
of simple roots. We take the roots of $B$ to be positive by
convention. Let $\Delta^{\vee}\subset R^{\vee}$ be the basis for 
the corresponding dual root system. Then we can define the Bruhat
ordering on ${\mathcal W}$. The longest element with respect to this
ordering of ${\mathcal W}$ is denoted by $w_0$.  A reductive group is
classified by these {\it root-data}, namely the character group,
1-parameter subgroups, the root system, co-roots and the ${\mathcal
W}$-invariant pairing.

 There is also a mapping from $X(T) \times Y(T) \to {\bz}, (\chi, \lambda) \mapsto \langle \chi, \lambda \rangle $ which is a dual pairing over $\bz$. The inner product on $X(T)$ gives one on $Y(T)$, say $(\lambda, \lambda')$. For $\lambda \in Y(T)$, define $\chi_{\lambda} \in X(T)$ by $\chi_{\lambda}(\lambda') = (\lambda, \lambda')$. Thus, $(\lambda, \lambda') = (\chi_\lambda, \chi_{\lambda'}) $. Since any two maximal tori are conjugate, for any 1-PS $\lambda$ of $G$, we have a well-defined norm $\parallel \lambda \parallel$ with $\parallel \lambda \parallel^2 \in {\bq}$.
 
 Following Hesselink, define $Y(G) = Hom({\mathbb G}_m, G)$ and define $Y(G,k')$ to be the one parameter subgroups of $G$ defined over $k'$ for any field $k' \supset k$. In particular, $Y(G) = Y(G,k)$. 
 
 Define $q(\lambda) := \parallel \lambda \parallel^2$ which defines a map $q: Y(T) \otimes \bq \to \bq$. The map $q$ extends to a $G$--invariant map from $M(G) \to \bq$, where $M(G) = (Y(G) \times {\bn}) /\sim$, where  $(\lambda,l) \sim (\mu,m)$ if $\lambda(t^m) = \mu(t^l)$. Note that $M(T) = Y(T) \otimes \bq$.

\bdefe\label{plambda} Let $\lambda \in Y(T)$. Define the associated parabolic subgroup $P_{_\lambda}$ as follows:
\beqa
P_{_\lambda} := \{g \in G \mid \lim_{t \to 0} \lambda(t) \cdot g \cdot\lambda(t)^{-1} ~ exists ~in~G\}
\eeqa
\edefe

 \bsem\label{kempf1ps} {\bf The Kempf-Rousseau 1-PS.} Let $V$ be a finite dimensional $G$--module and view it as a $T$-- module. Then we get a decomposition, $V = \oplus V_\chi$, summed over all characters $\chi \in X(T)$ such that $V_\chi \neq 0$. For elements $v \in V$, we express it as $v = \oplus v_\chi$. Define the {\em state of} $v$, $S_T(v) = \{\chi \mid v_\chi \neq 0 \}$.

For $\lambda \in Y(T) = Hom({\mathbb G}_m, T)$, we have $V = \oplus V_i$, where 
\beqa
V_i = \{v \in V \mid \lambda(a).v = a^i\cdot v ~~ \forall a \in k^* \}
\eeqa Thus, $V_i = \oplus V_\chi$, where the sum is over all characters $\chi$ such that $\langle \chi, \lambda \rangle = i$.

Let $V^q = \oplus_{i \geq q}V_i$. Then. $V^{q+1} \subset V^q$ and each $V^q$ is invariant under $P_{_\lambda}$ (see Definition \ref{plambda}). Thus, the group $P_{_\lambda}$ acts on the quotient $V^q/V^{q+1}$ and the quotient map 
\beqa
\pi: V^q \to {V^q \over V^{q+1}} 
\eeqa 
is $P_{_\lambda}$--invariant. Further, the unipotent radical $U(\lambda)$ acts trivially on $V^q/V^{q+1}$. 

Let $v \in V$. For a 1-PS $\lambda$ whose image is in some maximal torus $T'$, define
\beqa\label{v_o}
m(v, \lambda) := inf \{\langle\chi,\lambda \rangle \mid \chi \in S_{T'}(v) \}
\eeqa
Following Hesselink, we call  $m(v,\lambda)$  the ``measure of instability" (cf. \cite[2.2, page 77]{hesslink} and \cite[12.1]{kirwan}). The numerical criterion of stability is that $v \in V$ is {\em semistable} if and only if for every $\lambda$ we have $m(v.\lambda) \leq 0$. (Note that the $m(v,\lambda)$ here is different from the $\mu(v,\lambda)$ in Mumford's GIT and unfortunately called $\mu$ in \cite{rr}).

Fix $v_o \in V$ such that $0 \in {\overline{G.v_o}}$, i.e an {\em unstable} point.
G. Kempf (cf. \cite{kempf}) then showed that the function $\lambda \mapsto m(v_o, \lambda) / {\parallel \lambda \parallel}$ attains a maximum value on $Y(G)$ (see below the paragraph after Definition \ref{after}).

We may assume that this $\lambda \in Y(T)$. Then, $\lambda$ is uniquely determined among the indivisible 1-PS's of $T$ and is called the instability 1-PS for $v_o$. Put 
\beqa\label{j}
j := m(v_o, \lambda)
\eeqa
then $j = max \{q \mid v_o \in V^q \}$.\esem

Observe that the 1PS $\lambda$ determines a character $\chi_{_\lambda}$ on $P_{_\lambda}$ which is determined up to raising to a positive power. Further, if $p \in P_{_\lambda}$, then $P_{_{p \lambda^n p^{-1}}} = P_{_\lambda}$ for all $n > 0$. Moreover the associated characters satisfy the following relation: 
\beqa\label{q+ray}
\chi_{_{p \lambda^n p^{-1}}}^{r} = \chi_{_\lambda}^s
\eeqa
for some $r,s >0$ (similar to the definition of $M(G)$).

\bprop\label{rrkey} (\cite[Proposition 1.12]{rr}) Let $\lambda$ be the Kempf-Rousseau 1-PS for $v_o$. Then there exists a positive integer $r$ and a character $\theta_{\lambda}$ of $P_{_\lambda}$ determined up to equivalence as in \eqref{q+ray} by $v_o$ such that ${\overline{\pi({v_o}}}) \in {\bp}(V^j/V^{j+1})$ is {\em semistable} for the natural action of $P_{_\lambda}/U(\lambda)$ with the linearisation given by $\mathcal O(r) \otimes \mathcal O_{\theta_{\lambda}^{-1}}$. This gives a non-constant homogeneous function $f$ on $V^j/V^{j+1}$ such that $f(\pi(v_o)) \neq 0$ and $f(g.\pi(v)) = \theta_{\lambda}(g)^r f(\pi(v))$ for all $v \in V^j$ and $g \in P_{_\lambda}$. \eprop

\brem We may assume that $\theta_{\lambda}$ is a dominant character of $P_{_\lambda}$ (see \cite[Remark 1.13 ]{rr}).
\erem

\bsem {\bf Instability of points on a projective variety.} Let $S \hra {\bp}(V)$ be a projective $k$--variety with a linear action of $G$.  
Let $x \in S$ and let $x^{\ast}$ be a point above $x$ in the cone.   The value of $m(x^{\ast},\lambda)$ depends only on the point $x \in S$; therefore we simply write it  as 
$m(x,\lambda)$. We therefore have the following definition for points of a projective variety $S$:
\bdefe\label{after} A point $x \in S$ is unstable if $m(x,\lambda) > 0$ for some $\lambda \in Y(T)$. \edefe

For any $x \in S$, $x$ is unstable if and only if $q_G^{-1}(x) < \infty$, where 
\[
q_G^{-1}(x) := inf \{q(\lambda) \mid \lambda \in M(G), m(x,\lambda) \geq 1 \}
\]
and if
\[
\Lambda_G(x) := \{\lambda \in M(G) \mid m(x,\lambda) \geq 1, q(\lambda) = q_G^{-1}(x) \}
\]
then equivalently, $x$ is unstable if and only if  $\Lambda_G(x) \neq \emptyset$. 

The $T$ action on the vector space $V$ breaks it up into weight spaces indexed by characters $\epsilon_0, \ldots, \epsilon_n$ say. These $\epsilon_i$ which are a priori in $M(T)^*$ can be identified with elements in $M(T)$ by using $\parallel~\parallel$ whose square is $q$. Fix $x = (x_0:, \ldots, :x_n) \in V$ and let $\lambda$ be the {\em point closest to $0$} for the norm $q$ on the convex hull $C(x)$ of the set $\{\epsilon_i \mid x_i \neq 0 \}$ in $M(T)$. This closest point exists by the convexity and compactness of the closure of $C(x)$ in $M(T) \otimes {\mathbb R}$ and lies in fact in $C(x)$. Thus, the minimal is attained in the definition of $\Lambda_G$ (cf. \cite[Lemma 3.2]{hesslink}, or \cite[Lemma 1.1]{rr}). 

We recall the following lemma from \cite[Lemma 12.6]{kirwan}.
\esem

\blem If $\lambda \neq 0$ then $\Lambda_T(x) = \{\lambda/q(\lambda) \}$. \elem
Hesselink (cf. \cite[12.9]{kirwan}) defines the notion of ``optimality" of subgroups and one can always assume that the maximal torus $T$ is optimal for $x$. Recall the following result from \cite{kirwan}.

\bth\label{kr} (Kempf-Rousseau) For each unstable $x \in S$, there exists a unique parabolic $P(x)$ of $G$ such that $P(x) = P_{_\lambda} ~~\forall~ \lambda \in \Lambda_G(x).$ Further, if $T$ is optimal for $x$, and $\Lambda_T(x) = \{\lambda/q(\lambda)\}$, then $P(x) = P_{_\lambda}$.\eeth

\bdefe\label{thesetB} Consider the convex hull $C(F)$ in  $M(T)$ of a finite subset $F \subset \{\epsilon_0, \ldots, \epsilon_n \}$. An element in $C(F)$ is called  a {\em minimal combination of weights} if it is the closest point to $0$ in $C(F)$. Let ${\bf B}$ be the set of all minimal combination of weights in some positive Weyl chamber, for various $F$'s. \edefe

By the work of Kirwan (\cite[Theorem 12.26, Page 157]{kirwan}), we know that the variety $S$ has a natural stratification indexed by partially ordered set $\bf B$ defined above. For each $\beta \in \bf B$, we have a locally closed subvariety $S_{\beta}$ of $M$. The subvarieties $S_{\beta}$ are all $G$--invariant, i.e $G S_{\beta} = S_{\beta}$. Furthermore, 
\[
S = \bigsqcup_{\beta \in {\bf B}}  S_{\beta} 
\]
expressing $S$ as a disjoint union of the strata $S_{\beta}$.
The strata have a more precise geometric description: $S_0 = S^{ss}$ (i.e the $G$--semistable points), while if $\beta \neq 0$
\[
S_{\beta} = G Y^{ss}_{\beta}
\]
where $Y^{ss}_{\beta} = \{x \in S \mid \beta/q(\beta) \in \Lambda_G(x) \}$. If the variety $S$ is a nonsingular projective variety then the strata $\{S_{\beta} \mid \beta \in {\bf B}\}$ are all nonsingular and each stratum $S_{\beta}$ is in fact:
\[
S_{\beta} = G  \times ^{P_{_\beta}} Y^{ss}_{\beta}
\]

\section {The Bogomolov model}  {\em In this and the next section we will assume that $char(k) = 0$. We will also assume that $G$ is a semisimple connected algebraic group}. Notations are as in Section 5 above. Although the results from the later sections in this paper proves the theorem on semistability of the tensor product of two semistable Hitchin pairs over fields of characteristic zero as well as positive characteristics, we give a different argument in next two sections using the ideas of Bogomolov; this approach has the distinct advantage of making the ideas more transparent and in our opinion more geometric.

\bsem\label{section5} {\bf Towards Bogomolov's theorem.} Let $\chi \in X(T)$. Then it is well-known that there is a canonical parabolic subgroup $P(\chi)$ associated to $\chi$ and $\chi$ acts as a character on $P(\chi)$. It is the parabolic subgroup of $G$ generated by the maximal torus $T$ and the {\em root groups} $U_{_r}$ for $r \in R$ such that $(r, \chi) \geq 0$. In fact, if $\lambda$ is the 1 PS dual to $\chi$, then $P(\chi) \simeq P_{_\lambda}$ (see Definition \ref{plambda}). Let $L_{_-\chi}$ be the associated line bundle on $G/P(\chi)$. Let $W(\chi) \simeq H^0(G/P(\chi), L_{_-\chi})^*$. Then it is well-known that since we are in char $0$, $W(\chi)$ is an irreducible $G$--module with a non-zero vector $w_{\chi} \in W(\chi)$ (called a {\em highest weight vector}) unique up to constant multiples, such that for every $p \in P(\chi)$, one has:
\[
p.w_{\chi} = \chi(p).w_{\chi}
\]
Moreover, two such irreducible modules $W(\chi)$ and $W(\chi')$ are isomorphic as $G$--modules if and only if there is a $t \in \mathcal W$, such that $t.\chi = \chi'$. Let the notation $\bw(\chi)$ be as in \eqref{mathbbv}.

\bdefe\label{bogomodel} For any character $\chi \in X(T)$, $\chi \neq 0$ we define the Bogomolov model $A_\chi$ corresponding to $\chi$ to be the closure $\overline{G.w_{\chi}}$ in $\bw(\chi)$. \edefe 

We now return to the setting of Proposition \ref{rrkey}. We stick to the notations in \eqref{kempf1ps}. $V$ will be a finite dimensional $G$--module and $v_o \in V$ such that $0 \in {\overline{G.v_o}}$, i.e an {\em unstable} point.

Consider the finite dimensional irreducible $G$--module $W(r \theta_{\lambda})$ having highest weight vector $$w_{r \theta_{\lambda}} = w$$ with $T$--weight $r \theta_{\lambda}$. 

It is known that the stabilizer $\wp_{_\lambda}$ of the line  $k.w \subset \bw(r \theta_{\lambda})$  contains $P_{_\lambda}$ for the action on ${\bp}(W(r \theta_{\lambda}))$ and in fact $\{\wp_{_{\lambda}}\}_{_{red}} = P_{_\lambda}$. Since we work in characteristic $0$, we have $\wp_{_\lambda} = P_{_\lambda}$. 

Let $j$ be as in \eqref{j} and $V^j$ be as in Proposition \ref{rrkey}.  We first deduce a $P_\lambda$--morphism:
\beqa
\hat f : V^j \xrightarrow{\pi} V^j/V^{j+1} \xrightarrow{f} k
\eeqa
where $P_\lambda$ acts on $k$ by the weight $r \theta_\lambda$. 

Now consider the mapping $\psi_f: G \times V^j \to \bw(r \theta_{\lambda})$ given by 
\beqa\label{psi_f}
\psi_f(g, v) = g( \hat{f}(v) \cdot \{\hat{f}(v_o)^{-1}\cdot w\} )
\eeqa
This makes sense since ${\hat{f}}(v_o) \neq 0$, by Proposition \ref{rrkey}.

The group $P_{_\lambda}$ acts on $G \times V^j$ by $h (g,v) = (g h^{-1},h.v)$. Let 
\[
G \times^{P_{_\lambda}} V^j
\] be the quotient of $G \times V^j$ with respect to this action.
Let $G$ act on $G \times V^j$ by $x (g,v) = (xg,v)$. Now take 
\beqa\label{chiandlambda}
\chi = r.\theta_{\lambda}
\eeqa
Then one can check without much difficulty that $\psi_f$ is $P_{_\lambda}$--equivariant and we have a $G$--morphism:
\beqa
\psi_f : G \times^{P_{_\lambda}} V^j \to G \times^{P_{_\lambda}} k \to A_{\chi} \subset \bw(r \theta_{\lambda})
\eeqa
Also, $\psi_f(1,v_o) = w$.

Consider the map $G \times V^j \to V$, given by $(g,v) \to g.v$, which is $G$--equivariant. Let $Z$ be the image and $v_o \in Z$. Then $Z$ is {\em closed} and can be seen as follows (cf. \cite[Lemma 2.5]{rou}): 

Observe that $Z = G.V^j$ and $P_{_\lambda}.V^j \subset V^j$. Let $\alpha:G \times V \to G \times V$ be the map $(g,v) \to (g, g.v)$, and consider the following chain of maps:
\beqa
G \times V \xrightarrow{\alpha} G \times V \xrightarrow{\eta} {G / P_{_\lambda}} \times V \xrightarrow{pr_2} V
\eeqa 
with $\eta(g,v) = (g.P_{_\lambda},v)$. 

Then $G.V^j$ is the image by the composite $pr_2 .\eta. \alpha$ of $G \times V^j$. The map $\alpha$ is an isomorphism and hence $\alpha(G \times V^j)$ is closed in $G \times V$. Note that $P_{_\lambda}.V^j \subset V^j$. Therefore, we have $\eta^{-1}(\eta.\alpha(G \times V^j)) = \{(g,v) \mid g^{-1}v \in V^j \}$. This is closed in $G \times V$ and is isomorphic to $G \times V^j$ via the isomorphism $\alpha$. In other words, $\alpha(G \times V^j)$ is saturated for $\eta$. Since $\eta$ is a quotient morphism,  it follows that $\eta .\alpha(G \times V^j)$ is closed in ${G / P_{_\lambda}} \times V$. Therefore, since $G/P$ is proper, the image $G.V^j$ is closed in $V$ (cf. \cite[Lemma 2.5]{rou}).

The theorem of Bogomolov (cf. \cite[Theorem 2.7]{rou} and \cite[Page 287]{rou}) states that the map $\psi_f$ \eqref{psi_f} factors through $G.V^j$ by a $G$--morphism 
\[
\phi_f:G.V^j \to A_{\chi}
\] 

In fact, if $Y \subset I_G(V)$ is a closed irreducible $G$--subvariety then there exists a $\theta_{\lambda} \in X_{\mathbb Q}$ such that $Y \subset G.V^j$, where $j$ and $\theta_{\lambda}$ are as in Proposition \ref{rrkey}.

We now summarize the above discussion in the following key theorem.\esem

\bth\label{maintheoremofbogo}(Bogomolov) Let $I_G(V)$ be the subset of ${\mathbb V}$ consisting of the unstable points for the $G$ action. Let $Y$ ($Y \neq 0$) be a $G$--invariant closed subvariety of $I_G(V)$. Then there exists a $\chi \neq 0$, $\chi \in X(T)$ and a non-trivial $G$--morphism from $Y \to A_{\chi}$. 

More generally, if $Y = I_G(V)$, then there exists a filtration 
\beqa\label{bogofilt}
Y_0 = Y \supset Y_1 \supset \ldots \supset Y_m = \{0\}
\eeqa
of $Y$ by closed $G$--stable subvarieties such that for $i < m$, $Y_{i+1}$ is the intersection of the inverse images of $0$ by all the $G$--morphisms of $Y_i$ to models $A_{\chi_i}$.
\eeth

\noindent
\noindent {\it Proof}:
 See \cite[Theorem 2.7, page 284]{rou} for details.

\section{The main theorem in char $0$} 
The aim of this section as well as the main strategy in the arguments in the later sections is to use the interplay of Bogomolov instability of certain Higgs sections of associated bundles coming from a $G$--Hitchin pair and the Higgs instability of the $G$--Hitchin pair. This was indeed the strategy of Bogomolov and also Ramanan and Ramanathan. Let $G$ be a connected semisimple algebraic group.  

Let $V$ be a finite dimensional $G$--module and let $W = W(\chi)$ be as in \eqref{section5}. Let $A_{\chi}$ be as in Definition \ref{bogomodel}. 
 
\bsem{\bf Reduction to the Kempf-Rousseau parabolic.} Denote by ${\bp}(A_{\chi})$ the image in ${\bp}(W)$ of $A_{\chi} - (0)$; it is a closed subset of the projective space and therefore a projective variety. The group $G$ acts on ${\bp}(W)$ and ${\bp}(A_{\chi})$ is an orbit for this action.\esem

\bth\label{maintheorem1} Let $(E,\theta)$ be a principal $G$--Hitchin pair. Let $\sigma$ be a Bogomolov unstable Higgs--section of $E({\mathbb V})$ for the induced Higgs structure $\theta_V$. Then,
\begin{enumerate}
\item There exists a $\chi \in X(T), \chi\neq 0$ (which we may assume to be dominant), and a non-zero geometric Higgs section $s$ of the associated Hitchin
 scheme
 $E(A_\chi)$.
\item The projected section $s_1$ of $E({\bp}(A_\chi))$ on a non-empty open $U \subset X$ extends to $X$.

\end{enumerate}
\eeth
\noindent
\noindent {\it Proof}:
 
\noindent
{(1)}: Let $Y = I_G(V) \subset {\mathbb V}$ be the $G$--subvariety of unstable points of ${\mathbb V}$. By the main theorem of Bogomolov (see \cite{bogo} and \cite[Corollaire 2.8]{rou}), there exists a filtration \eqref{bogofilt} above. Taking the corresponding associated Hitchin schemes we have a filtration:
\beqa\label{largestindex}
~~~~~E(Y) = E(Y_0) \supset E(Y_1) \supset \ldots \supset E(Y_m) = X \times \{0\}
\eeqa
{\em Let $i$ be the largest index such that $\sigma(X) \subset E(Y_i)$}. Then by the canonical property of the Bogomolov model and the instability of the section, we have 
\begin{itemize}
\item a $G$--morphism $\varphi:Y_i \to A_{\chi} \subset {\mathbb W}$ for a suitable character $\chi = \chi(i)$ and 
\item furthermore, the induced map $E(\varphi): E(Y_i) \to E(A_{\chi}) \subset E({\mathbb W})$ when evaluated on the subset $\sigma(X) \subset E(Y_i)$, has the property that $ E(\varphi)(\sigma(X)) \neq 0$.
\end{itemize}
 Define
\beqa
s := E(\varphi)(\sigma)
\eeqa
Since the section $\sigma: X \to E({\mathbb V})$ is a Higgs section, and since the image of $\sigma$ lies in $Y_i$, the section $\sigma:X \to E(Y_i)$ is a geometric  Higgs section of the associated Hitchin
 scheme
 $E(Y_i)$. Further, by Proposition \ref{assocmorphs} the map $E(\varphi)$ is a map of Hitchin
 schemes
 and the induced map $s = E(\varphi)(\sigma)$ is therefore a non-zero geometric Higgs section of $E(A_{\chi}) \subset E({\mathbb W})$.

This proves (1).

\noindent
\linebreak
{(2)}: Let $U = \{x \in X \mid s(x) \neq 0 \}$. By projecting to ${\mathbb P}(A_{\chi})$ we get a section of $E({\mathbb P}(A_{\chi}))$ on the open subset $U$. Since ${\mathbb P}(A_{\chi})$ is projective, the section uniquely extends to  a section of $E({\mathbb P}(A_{\chi}))$ on $X$. \begin{flushright} {\it QED} \end{flushright}

\bth\label{theparabolichiggsred}
The construction of the extension $s_1:X \to E({\mathbb P}(A_{\chi}))$ in Theorem \ref{maintheorem1} gives rise to a  $P_{_\lambda}$--reduction of the principal Hitchin pair $(E, \theta)$, to say, $E_{_{P_{_\lambda}}} \subset E$. Furthermore,  the reduction to $P_{_\lambda}$ is compatible with the Higgs structure on $(E, \theta)$. 
\eeth 
\noindent 
\noindent {\it Proof}:
 By the construction of the Bogomolov model and the discussion in \eqref{section5}, we have a dominant character $\chi$, a $1$-PS $\lambda$ and the character $\chi$ and $\lambda$ are related by \eqref{chiandlambda}. Since ${\bp}(A_{\chi}) = G/P_{_\lambda}$, the section obtained above gives a reduction of structure group of $E$ to the parabolic subgroup $P_{_\lambda}$. To get the Higgs structure on $E_{_{P_{_\lambda}}}$ we proceed as follows (see Remark \ref{generichiggsred} and Remark \ref{caution}). 

Let $H = Stab_G(v)$ for the $G$--action on $A_{\chi}$. Observe that $H \subset P_{_\lambda} \subset G$. The scheme $G/H \subset A_{\chi}$ is a subscheme. Now the section $s:X \to E(A_{\chi})$ has the property that the image $s(\xi)$ of the generic point $\xi \in X$,
lies in $E(G/H)$. 

If we base change to $Spec(K) \subset X$ the $G$--equivariant Higgs morphism $\phi: E \to A_{\chi} \times X$ we are in the setting of Lemma \ref{fibreproduct} since inverse images can be realized as fibre products. We then immediately obtain the consequence that the $K$--subscheme $\phi^{-1}(eH \times K) \subset E_K$ is a Higgs subscheme (see Remark \ref{generichiggsred}). By standard arguments, we see that  $\phi^{-1}(eH \times K) = (E_{_H})_K$ is a $H$--reduction of $E_K$.

Thus, the section $s$ provides a Higgs $H$--reduction over $K$, and hence a Higgs $P_{_\lambda}$--reduction over $K$.

As observed earlier (proof of Theorem \ref{maintheorem1}), the underlying reduction of structure group $E_{_{P_{_\lambda}}}$ of $E$ to $P_{_\lambda}$ extends as a reduction of structure group to the whole of $X$. Hence by Lemma \ref{genhiggs}, it follows that the reduction of structure group to $P_{_\lambda}$ is a Higgs reduction on the whole of $X$. \begin{flushright} {\it QED} \end{flushright}

\bsem{\bf Higgs semistability and associated bundles.}
Let $M$ be a finite dimensional $G$--module such that $\rho:G \to SL(M)$ is a representation. Let $Q \subset SL(M)$ be a maximal parabolic and let $L_{_\gamma}$ be an ample line bundle coming from a dominant character $\gamma$. Let $SL(M)/Q \subset {\mathbb P}(W)$ be the embedding defined by $L_{_\gamma}$. Let $s: X \to E(SL(M)/Q)$ be a section and suppose that $s$ lifts to a section of the associated bundle on the cone of $SL(M)/Q$, say $E(\hat{C})$. Suppose that $s$ is an unstable section of $E(W)$. Thus, 
\[
s(X) \subset \hat {C} \cap I_G(W)
\]
hence by Theorem \ref{maintheorem1}, $s(X) \subset Y_i \cap \hat {C}$ (the index $i$ as defined after \eqref{largestindex}). 

Therefore, the induced $G$--morphism $Y_i \to A_{\chi}$ for some $\chi$,  is such that the image of the general point of $\xi \in X$ maps to the highest weight vector $w_{\chi}$ in the model $A_{\chi}$ and the stabilizer of the line $k w_{\chi}$ in ${\mathbb P}(A_{\chi})$ is the Kempf-Rousseau parabolic $P_{_\lambda}$ (where $\chi$ is related to $\lambda$ as in \eqref{chiandlambda}).

\bth\label{ramram} (See Theorem \ref{maintheorem} for the result in positive characteristics as well) Let $(E,\theta)$ be a Higgs principal $G$--Hitchin pair, with $G$ semisimple. Suppose that $\rho: G \to GL(M)$ be a representation. If $(E(M), \theta_M)$ is a Higgs {\em unstable} $GL(M)$--Hitchin pair  of degree $0$ so is $(E, \theta)$.\eeth

\noindent
\noindent {\it Proof}: We begin with the $GL(M)$--Hitchin pair $(E(M), \theta_M)$ which is {\em Higgs unstable}. By Definition \ref{stabilitynotions}, this implies that there is a maximal parabolic $Q \subset GL(M)$ and a dominant character $\eta$ on $Q$ together with a Higgs reduction $s:X \to E({GL(M)\over Q})$ such that  the pull-back by $s$ of the associated  line bundle $L_{\eta}$ has $deg(s^*(L_{\eta})) > 0$. 

Choose $m \gg 0$ such that 
\beqa\label{genus}
deg(s^*(L_{m\eta})) > {\tt g}
\eeqa 
where {\tt g} = genus(X).

Consider the {\em dual} $L_{\eta}^{\vee}$. By the convention (cf. Definition \ref{stabilitynotions}), $L_{\eta}^{\vee}$ is ample. Let $V = H^0(L_{m\eta}^{\vee})^*$. Then we have the Pl\"ucker embedding ${GL(M)\over Q} \subset {\mathbb P}(V)$. Observe that if $Q$ fixes the subspace $M_1 \subset M$, and $dim(M_1) = r$, then, $V \subset Sym^{m}(\wedge^{r}(M)))$. 

Write ${\mathbb P}(V) \simeq GL(V)/P{_{_\ell}}$, and since $GL(M)/Q \subset GL(V)/P{_{_\ell}}$ and the subgroup $Q$ fixes the line $\ell \subset V$. Hence $Q \subset P{_{_\ell}}$.
The Higgs section $s:X \to E(({GL(M)\over Q})$ gives a $Q$-Hitchin pair structure on $E_{_Q}$. Hence  by composing with the inclusion $E({GL(M)\over Q}) \hookrightarrow E({\mathbb P}(V))$, we get, via extension of structure group by the inclusion $Q \subset P{_{_\ell}}$, a $P{_{_\ell}}$--Hitchin pair structure on $E_{_Q} \times ^{Q} P{_{_\ell}} \simeq E_{P{_{_\ell}}}$.

By the discussion in  Example \ref{eg}, this section therefore gives a Higgs subbundle $L' \subset E(V)$. Furthemore, by \eqref{genus}, $deg(L') > {\tt g}$ which in particular implies that $L' = \co_X(D)$ for an effective divisor $D$. 

In other words, $s$ induces a section $s:\co_X(D) \hookrightarrow E(V)$ (we use the same notation for $s$) which implies by the definition of Bogomolov instability that $s$ is a Bogomolov unstable Higgs section of $(E(V), \theta_V)$. We may view the Higgs section $s$ as a section of the geometric Hitchin scheme $(E(\mathbb V), \theta_V)$ which is such that the zeroes of $s$ coincide with the effective divisor $D$. 

Since $s$ is a Bogomolov unstable Higgs section, we use Theorem \ref{maintheorem1} to get a Bogomolov model $A_{\chi}$ together with a Higgs section of the Hitchin scheme $E(A_{\chi})$. That is, the geometric section $s:X \to E(\mathbb V)$ factors through $s:X \to E(A_{\chi})$. The section $s$ maps $X - D$ to $E(A_{\chi} \smallsetminus \{0\})$. 

By Theorem \ref{theparabolichiggsred}, $s$ extends to a section \[t:X \to E(\bp(A_{\chi})) = E\big({G/ P_{_\lambda}}\big)\] which moreover gives a Higgs reduction of structure group to the Kempf-Rousseau parabolic $P_{_\lambda}$. 

The character $\chi$ and $\lambda$ are related by \eqref{chiandlambda}. The dominant character $\chi$ is a character of $P_{_\lambda}$ and gives rise to a line bundle $L_{\chi}$ on $G/P_{_\lambda}$ such that $L_{\chi}^{\vee}$ is {\em ample}.

Note that $t$ is obtained by composing $E(A_{\chi} \smallsetminus \{0\}) \to E(\bp(A_{\chi})) = E\big({G/ P_{_\lambda}}\big)
$ with $s$ and the geometric bundle underlying $t^*(L_{\chi}) $ is gets identified with
\beqa
t^*\big(E(A_{\chi} \smallsetminus \{0\})\big)
\eeqa
Hence the induced section $t$ in fact imbeds $\co_X(D) \subset t^{*}(L_{\chi})$. This implies that $deg(t^{*}(L_{\chi})) > 0$. By Definition \ref{stabilitynotions}, this gives the Higgs instability of $(E, \theta)$. \begin{flushright} {\it QED} \end{flushright}\esem

\section{Theorems in positive characteristics}
Let the ground field $k$ be algebraically closed of arbitrary characteristics. The notations are as in Section 5. Let {\em $G$ be a connected semisimple algebraic group}. Let $K = k(X)$ be the function field of the base curve and let $\xi \in X(K)$ be the generic point of $X$. Let $S \hra {\bp}(V)$ be a projective $k$--variety with a linear action of $G$ and where $V$ is a low height $G$--module (see \eqref{lowheight}). We observe that inputs from (\cite{rr}) and (\cite{imp}) (see Theorem \ref{imp} below) allows one to conclude that the Kempf-Rousseau parabolic is defined over $K$. The existence of a reduction of structure group to $P_{_\lambda}$ is concluded rather scheme--theoretically in \cite{rr} and does not adequately reflect the geometry. This makes it almost impossible to generalize the strategy to the setting of Hitchin pairs.

The Kirwan stratification also gives a geometric description of the strata and we derive a geometric realization of the Kempf reduction in the case of Hitchin pairs. We remark that even without the Higgs structures, the proofs that we give here makes the entire theory more transparent. We have therefore taken the opportunity to briefly expound the central point in the proof of the main theorem of \cite{imp}.

The following ideas underlie the proof in positive characteristics. 

\begin{enumerate}
\item[\bf(a)] Let $S \hra {\bp}(V)$ be a projective $k$--variety with a linear action of $G$ and where $V$ is a low height $G$--module. Let $(E, \theta)$ be a principal $G$--Hitchin pair. Let $s$ be an GIT unstable section of the underlying associated fibration
 $E(S)$ (see Definition \ref{gitunstable} below). For an unstable $K$--point $s(\xi)$ the rationality of the Kempf-Rousseau parabolic follows from Theorem \ref{imp} (cf. \cite{imp} and \cite{rr}).

\item[\bf(b)] Once this rationality is achieved, then Hesselink\cite[Theorem 5.5, page 82]{hesslink} shows that the Kirwan strata $S_{\lambda}$ containing $s(\xi)$ is actually defined over $K$; the proof of this uses the data from (a).

\item[\bf(c)] From this and Kirwan's description of the strata we get the fact that the 
morphism $S_{\lambda} \to G/P_{_\lambda}$ is defined over $K$.

\item[\bf(d)] This gives the reduction of structure group to the Kempf-Rousseau parabolic $P_{_\lambda}$.

\item[\bf(e)] The Higgs geometry defined earlier then uses this $K$--morphism
to get a geometric description of the Higgs reduction to the Kempf-Rousseau parabolic $P_{_\lambda}$.

\item[\bf(f)]  Degree computation for this ``Higgs" reduction  to the Kempf-Rousseau parabolic then follows easily and gives the required semistability results.
 
\end{enumerate}

\bsem\label{lowheight}{\bf Low height representations.} We recall (cf. \cite[Page 7]{bapa}) the definition and some salient properties of  a low height representation $\rho:G \to SL(V)$. Recall the notations from Section 5. 

Observe that $V$ can be written as
direct sum of eigenspaces for $T$. On each eigenspace the torus $T$ acts by a
character. These are called the {\em weights} of the representation. A weight $\lambda$ is called {\em dominant} if
$(\lambda,\alpha_i^{\vee})\geq 0$ for all simple roots $\alpha_i\in
\Delta^+$. A weight $\lambda$ is said to be ``$\geq$" another weight $\mu$ if
the difference $\lambda-\mu$ is a non-negative integral linear combination
of simple roots, where the difference is taken with respect to the natural
abelian group structure of $X(T)$. The {\em fundamental weights} $\omega_i$
are uniquely defined by the criterion
$(\omega_i,\alpha_j^{\vee})=\delta_{ij}$.  The {\em
height}  of a root is defined to be the sum of the
coefficients in the expression $\alpha~=~\Sigma k_i\alpha_i$. We extend
this notion of height linearly to the weight space and denote this function
by $ht(~)$. Note that $ht$ is defined for all weights but need not be an
integer even for dominant weights. We extend this notion of height to
representations as follows:

\bdefe\label{height}(cf. \cite{imp}) 
\begin{enumerate}
\item Given a linear representation $V$ of $G$, we define the
{\em height} of the representation $ht_G(V)$ (also denoted by $ht(V)$ if
$G$ is understood in the given context) to be the maximum of
$2ht(\lambda)$, where $\lambda$ runs over dominant weights
occurring in $V$.  

\item A linear representation $V$ of $G$ is said to be a {\em low
    height} representation if $ht_G(V)<p$, and a weight $\lambda$ is of low height if $2ht(\lambda)<p$.  
\end{enumerate}
\edefe

Then we have the following theorem (cf. \cite{imp}, \cite{mour} and \cite{serre})

\bth
Let $V$ be a linear representation of $G$ of low height. Then
$V$ is semisimple. 
\eeth

\bcor Let $V$ be a low height representation of $G$ and $v \in V$ an
element such that the $G$-orbit of $v$ in $V$ is closed. Then $V$ is a
semisimple representation for the reduced stabiliser $G_{v,red}$ of $v$.
\ecor

\brem In Serre (\cite{serre} and \cite{mour}), the notation for $ht_G(V)$ is simply $n(V)$. \erem
\esem

\bsem{\bf Rationality issues.} Recall the following definitions from Hesselink (\cite{hesslink}) and the book by Kirwan (\cite{kirwan}). We observe also that the notable difference between the treatment in Hesselink (\cite{hesslink}) and Kempf (\cite{kempf}) is that \cite{hesslink} works over arbitrary fields while Kempf assumes that $K$ is perfect.\esem

Recall that $Y(G)$ is the set of 1PS's $\lambda: {\mathbb G}_m \to G$ of $G$ defined over $k$ (see Section 5).

\bdefe Let $V$ be a $G$--module. Let $\lambda \in Y(G)$. A point $v \in V(K)$ is called $\lambda$--unstable if $\lambda$ drives $v$ to zero. The point $v$ is called $L$--unstable for an extension $L/K$ if it is $\lambda$--unstable for a $\lambda \in Y(G,L)$, {\em where we denote by $Y(G,L)$ the subset of elements of $Y(G)$ which are defined over the field $L$}. 

A subset $\hat{S} \subset V(L)$ is called uniformly $L$--unstable if there exists a $\lambda \in Y(G,L)$ such that all $s \in \hat{S}(L)$ are $\lambda$--unstable. \edefe

We now summarize the main result of \cite{imp} in the following theorem:

\bth\label{imp} Let $\rho: G \to SL(V)$ be a low height representation and let $Q \subset SL(V)$ be a maximal parabolic subgroup and let us denote the homogeneous space $SL(V)/Q$ by $S$.

Let $\chi$ be a dominant character on $Q$ and let $L_{\chi}^{\vee}$ be the ample line bundle on $S$ embedding $S \hra {\bp}^r$. Let $K/k$ be an extension field. Let $m \in S(K)$ be a $\bar K$--unstable point for the $G$--action. Let $P(m)$ be the Kempf-Rousseau parabolic subgroup of $G$ given by the 1PS $\lambda$ which ``optimally drives m to zero" over $\bar {K}$, i.e $$P(m) = P_{_\lambda}$$ Then, $P(m)$ is defined over $K$. In other words, $m$ is a $K$--unstable point. \eeth
\noindent {\it Proof}:
 This theorem is  a generalization of a theorem in \cite[Theorem 2.3.]{rr} where similar rationality questions are addressed. In \cite{rr} the crucial assumption is that the action is {\em strongly separable}, which essentially ensures that the isotropy subgroups for the action of $G$ at any point are (absolutely) reduced.

The assumption on the height of the representation allows for this generalization. Since the proof for the most part   follows \cite{rr}, we will follow closely the notations and recall the relevant details from there. 

The assumption of separability of the action in \cite[Theorem 2.3]{rr} is used in \cite[Page 279, paragraph 1]{rr}. Instead, one uses the hypothesis of Theorem \ref{imp}, namely 
\begin{enumerate}
\item the low height property of $G \to SL(V)$.
\item the fact that the variety where the action is being studied is a {\em Grassmannian} and not an arbitrary projective variety.
\end{enumerate}
We then use this together with \cite[Lemma 3.2 and Proposition 3.3]{imp}. With these changes in place, the rest of the proof of \cite[Theorem 2.3]{rr} goes through without any difficulty. 

The key step in the proof, as in all rationality questions, is to show that $P(m)$ is defined over the separable closure $K_s$. Then the $\lambda$ will lie in $Y(G, K_s)$ and the point $m$ will firstly be shown to be $K_s$--unstable and then a Galois descent argument will show that it is actually $K$--unstable. \begin{flushright} {\it QED} \end{flushright}

Recall another {\em rationality} theorem from \cite[Theorem 5.5 (a)]{hesslink}:

\blem\label{rationality} Let $\hat{S} \subset V$ be a closed subset. Then $\hat{S}$ is uniformly $K_s$--unstable if and only if $\hat{S}$ is uniformly $K$--unstable. \elem

From these lemmas and the Theorem \ref{imp} we get the following key corollary. 

\bcor Let $\rho:G \to SL(V)$ be a low height representation and let $S = SL(V)/Q$ and  $m \in S(K)$ be a $\bar K$--unstable point driven to zero ``optimally" by $\lambda \in Y(G, \bar K)$. Let $S_{\lambda} \subset S$ be the strata containing $m$. Then the locally closed subvariety $S_{\lambda} \subset S$ is defined over $K$. \ecor

\noindent {\it Proof}:
 Observe that since $m \in S(K)$ is $\bar K$--unstable by the low height assumption on $V$ and by Theorem \ref{imp}, it follows that $m \in S(K)$ is $K$--unstable and the Kempf-Rousseau parabolic $P(m)$ and the $\lambda$ are both defined over $K$. That is $\lambda \in Y(G, K)$.

The strata $S_{\lambda}$ containing $m$  is a priori uniformly $\bar K$--unstable and since $m \in S_{\lambda}$, it follows by Theorem \ref{imp} that the 1-PS $\lambda$ actually lies in $Y(G,K_s)$.  Hence the Kempf-Rousseau parabolic $P(m) = P_{_\lambda}$ is also defined over $K_s$. That is, $S_{\lambda}$ is uniformly $K_s$--unstable. 

Since $P_{_\lambda}$ and $\lambda$ are both defined over $K$ itself, by \cite[Theorem 5.5]{hesslink} it follows that the strata $S_{\lambda}$ is uniformly $K$--unstable. In other words, $S_{\lambda}$ is defined over $K$.

By \cite[Proposition 6.1]{hesslink}, since $m \in S(K)$ is a $K$--unstable point  it implies that there is a maximal $K$--unstable subset $Y^{ss}_{\lambda} \subset S$  with $m \in Y^{ss}_{\lambda}$ such that all its points are driven to zero by $\lambda \in Y(G,K)$. By \cite{hesslink} again each $Y^{ss}_{\lambda}$ is invariant under the action of the Kempf-Rousseau parabolic $P_{_\lambda}$. 

As we have seen in our description of Kirwan's stratification, we know that $G \times ^{P_{_\lambda}}(Y^{ss}_{\lambda}) = S_{\lambda}$ is a locally closed subvariety of $S$ defined over $K$. In conclusion, we have shown that $S_\lambda$ is defined over $K$ and we have a $K$--morphism $S_\lambda \to G/P_{_\lambda}$ with fibres isomorphic to $Y^{ss}_\lambda$. \begin{flushright} {\it QED} \end{flushright}

\bsem{\bf GIT instability and unstable reductions.} 

\bdefe\label{gitunstable} Let $E$ be a principal $G$--bundle on $X$. Let $Y$ be an affine $G$--variety. Let $s:X \to E(Y)$ be a section of the associated fibration.  We say the section is {\em GIT unstable} if the evaluation $s(\xi)$ of the section at the generic point $\xi \in X(K)$ is a GIT unstable point in $E(Y)(K)$ for the $G_K$--action on $E(Y)_K$. If $Y$ is a projective variety on which there is a $G$--linearized action with respect to a very ample line bundle $\cl$, and $s:X \to E(Y)$ a section, we call it GIT unstable if $s(\xi)$ is a GIT unstable point in the cone over $E(Y)_K$.
\edefe 

\brem Observe that this definition makes perfect sense even if $X$ is not projective. \erem

\brem Let $X$ be projective and $Y$ be an affine $G$--variety. A section $s$ is GIT unstable in the sense of Definition \ref{gitunstable} if and only if  for some point $x \in X$ image $s(x) \in E(Y)_x$ is GIT unstable for the $G$--action. In other words, the section $s$ is Bogomolov unstable in the sense of Definition \ref{bogohiggs}. Of course, it must be noted that these notions make sense even when we have a Higgs structure. \erem

\blem\label{clarificationofgit} Let $E$ be a principal $G$--bundle. Let $\rho:G \to SL(W)$ be a linear representation. Consider the associated principal $SL(W)$--bundle $E(SL(W))$. Suppose that there exists a  maximal parabolic $Q \subset SL(W)$ and a dominant character $\eta$ of $Q$  together with a reduction of structure group $s:X \to E(SL(W)/Q)$  such that 
\beqa\label{degbound1}
deg(s^*(L_{\eta}))> 0
\eeqa 
Then, the section $s$ is a GIT unstable section in the sense of Definition \ref{gitunstable} for the ample line $L_{\eta}^{\vee}$ on $SL(W)/Q$. \elem 
{\it Proof}: This follows immediately from Lemma \ref{degandss}. 
\esem

\brem We note that in the Lemma \ref{clarificationofgit} the parabolic subgroup $Q$ gives a destabilizing vector subbundle for the vector bundle $E(W)$. \erem

We now prove the following theorem:

\bth\label{maintheorem} Let $(E, \theta)$ be a Higgs semistable principal $G$--Hitchin pair, $G$ being semisimple. Let $\rho:G \to SL(W)$ be a low height representation. Then the Hitchin pair $(E(W), \theta_W)$ is Higgs semistable.\eeth

\noindent {\it Proof}:
  The proof breaks up into two parts.

\noindent
{\it {Higgs compatible Kempf-Rousseau parabolic reduction}}:

Let $Q \subset SL(W)$ be a  maximal parabolic  and $\eta$  a dominant character of $Q$. Suppose that we are given a Higgs-reduction of structure group $s:X \to E(SL(W)/Q)$. There are two possibilities:
{\em either the section is GIT semistable or it is GIT unstable} in the sense of Definition \ref{gitunstable}.

In the first case, by Lemma \ref{degandss}, since $L_{\eta}^{\vee}$ is ample, 
\beqa\label{degbound}
deg(s^*(L_{\eta})) \leq 0
\eeqa
In the second case it is more subtle and we proceed as follows. 
Let $$S := SL(W)/Q$$ Let $\xi \in X$ be the generic point and let $s(\xi) \in E(S)(K)$ be GIT unstable for the action of $G$. Let us denote the image $s(\xi)$ by $m$.

 By the assumption of low height the 1-PS $\lambda$ and the Kempf-Rousseau parabolic $P_{_\lambda}$ are both defined over $K$ (by Theorem \ref{imp}). Hence, $m \in S_{\lambda}$ with $S_{\lambda}$ defined over $K$ (see the end of Section 5 and Theorem \ref{kr} for the notations). Let $E_\xi$ be the generic fibre of $E \to X$. Now view the restriction of Higgs section $s$ to $Spec(K)$:
\beqa
\xymatrix{
E_{\xi} \ar[dr]{} \ar[rr]^{s}& & S_{\lambda} \times K  \hra S \times K \ar[dl]^{}\\
& Spec(K) &
}
\eeqa 
\noindent
as a $G$--equivariant morphism of Hitchin schemes over $Spec(K)$. Composing with the canonical projection $S_\lambda \to G/P_{_\lambda}$, we get the morphism $t$ of Higgs $K$--schemes:
\beqa
\xymatrix{
E_{\xi} \ar[dr]{} \ar[rr]^{t}& & G/P_{_\lambda} \times K \ar[dl]^{} \\
& Spec(K) &
}
\eeqa 
Taking the inverse image of $t^{-1}(eP_{_\lambda} \times K)$, we get $E_{P,\xi}\subset E_{\xi}$ which is a $P_{_\lambda}$--subscheme of $E_{\xi}$. 

Observe further that since $X$ is a curve, the reduction to $P_{_\lambda}$ extends as a usual reduction (i.e without the Higgs structure). What we have shown is that $E$ has a $P_{_\lambda}$--reduction which is generically Higgs. By Lemma \ref{genhiggs}, it follows that this reduction is a global Higgs reduction.

\noindent
{\it {The degree computations}}:

As in the proof of Theorem \ref{ramram}, using the ampleness of the {\em dual} $L_{\eta}^{\vee}$ and taking sections for a suitable power $m$ of $L_{\eta}^{\vee}$ we get the Pl\"ucker embedding of $S \subset {\mathbb P}(V)$. 

The parabolic subgroup $P_{_\lambda}$ being defined over $K$ gives a Higgs reduction of
structure group of $(E,\theta)$ to a parabolic $P_{_\lambda}$ of $G$. Let $V = \bigoplus
V_i$ be the weight space decomposition of $V$ with respect to $\lambda$. Let $V^q = \oplus_{i \geq q}V_i$ (see the discussion in \eqref{kempf1ps}).

Let \[
j = \mu(m, \lambda) = min \{i \mid m ~~has~~ a~~ non-zero~~ component~~ in~~ V_i \} = max \{q \mid m \in V^{q}\}\]
Let $m = m_0 + m_1$, with $m_0$ of weight $j > 0$ and $m_1$ the sum of
terms of higher weights. In other words, in the projective space ${\mathbb
P}(V)$ we see that $\lambda(t) \cdot m \longrightarrow m_0$. It is not too hard to see
that we have an identification of the Kempf-Rousseau parabolic subgroups
associated to the points $m$ and $m_0$ which is therefore simply denoted by $P_{_\lambda}$. (cf. \cite[Proposition 1.9]{rr}).

In the generic fibre $E(V)_{\xi}$ we have the projection
\beqa\label{mtomnought}
\bigoplus_{i \geq j} V_i \longrightarrow V_j
\eeqa
which takes $m$ to $m_0$. This gives a line sub-bundle $L_0$ of $E(V_j)$  associated  to the point $m_0$ as well as a nonzero map $s^*(L_{\eta}) \to L_0$. 

Now by Proposition \ref{rrkey},  $m_0$ is
in fact GIT semistable for the action of $P_{_\lambda}/U$, i.e for the Levi quotient of $P_{_\lambda}$, for a
suitable choice of linearisation obtained by twisting the action by a
dominant character $\chi$ of $P_{_\lambda}$ (see \eqref{chiandlambda}).  Let $E_{_{P_{_\lambda}}}$ be the $P_{_\lambda}$--bundle obtained by the Higgs reduction to $P_{_\lambda}$.  Let $L_{\chi} = E_{_{P_{_\lambda}}}(\chi)$. Since $(E,\theta)$ is Higgs semistable by assumption, we have
\beqa\label{contra}
deg(L_{\chi}) \leq  0
\eeqa
The GIT semistability of the point $m_0$ for the action of $P_{_\lambda}/U$ with respect to this new linearisation,  forces, by Lemma \ref{degandss}, the following degree inequality (we skip the details which are essentially in the closing parts of the proof of \cite[Proposition 3.13]{rr}):
\beqa\label{contra0}
deg(L_0 \otimes L_{\chi}^{\vee}) \leq 0
\eeqa
By \eqref{contra} $deg(L_{\chi}^{\vee}) \geq 0$. Hence by \eqref{contra0}, $deg(L_0) \leq 0$. This implies $deg(s^*(L_{\eta})) \leq  0$. 
This shows that in either case, we get the inequality \eqref{degbound}. This proves  the Higgs semistability of $(E(SL(W),\theta_W)$.\begin{flushright} {\it QED} \end{flushright}

\brem We remark that if we follow the proof strategy of Theorem \ref{ramram} in the positive characteristic case, the prime bounds that are forced are much bigger than those imposed by {\em low height} considerations. \erem 

\bth Let $(V_1, \theta_1)$ and $(V_2, \theta_2)$ be two Higgs semistable Hitchin pairs with $det(V_i) \simeq \co_X, i = 1,2$. Suppose that the ground field $k$ has characteristic $p$ such that
\[
rank(V_1) + rank(V_2) < p + 2
\]
Then the tensor product $(V_1 \otimes V_2, \theta_1 \otimes 1 + 1 \otimes \theta_1)$ is also Higgs semistable. 
\eeth

{\it Proof}: This is immediate from height computations for the tensor product representations and Theorem \ref{maintheorem} above. See for example \cite[5.2.5]{serrebourbaki}. \begin{flushright} {\it QED} \end{flushright}

\section{Polystability of associated bundles}
The ground field $k$ has arbitrary characteristics in this section. Let $G$ be a connected reductive algebraic group.  Let $T$ be the maximal torus of $G$ and $W$ be a finite dimensional $G$-module.
Further, let ${X}(T)$ be the free abelian group of characters of $T$
and $\mathcal S$ be the set of distinct characters that occur in $W$.

For every subset $S \subset \mathcal S$ we have the following map:
\[
\nu_S: {\bf Z}^{|S|} \longrightarrow {X}(T)
\]
given by $e_s \longrightarrow \chi_s$. Let $g_S$ be the g.c.d of the maximal minors of the map $\nu_S$ written under the fixed basis. 

For any vector  $w\in W$, consider the
subset $S_w \subset \mathcal S$, consisting of characters that occur in $w$
with nonzero coefficients. i.e., if $w = \sum a_{\chi}(w)e_{\chi}$, then
\[
S_w = \{{\chi \in \mathcal S} | a_{\chi}(w) \neq 0 \}
\]
Then we recall the following:

\blem (cf. \cite[Lemma 6]{bapa}) The characteristic of the field, $p$ does not divide $g_{S_w}$ if 
and only if the action of $T$ on the vector $w$ is separable.  
\elem

Define
\[
{\sf p}_T(W) := \{{\rm largest~prime~which~divides}~g_S |\forall S
\subset \mathcal S\} 
\]
\bdefe\label{sepind} Let $\rho:G \longrightarrow SL(W)$ be a finite dimensional
representation of $G$. Define the {\it separability index}, $\psi_{_G}(W)$
of the representation as follows:
\[
\psi_G(\rho) = \psi_{_G}(W) = max\{ht_{_G}(W), ~{\sf p}_T(W)\}
\]
\edefe

The notion of separability index was first defined in (\cite{bapa}). The reader is referred to (\cite[Section 4]{bapa}) for the details.

\bdefe\label{lowsep} The module $W$ is said to be with {\em low separability index} if the characteristic of the ground field $k$ is either {\em zero} or $p$ which satisfies $p > \psi_{_G}(W)$. \edefe
We recall the results proved in \cite[Proposition 5, p 16]{bapa} along with a key result from \cite{bardsley}.

\bprop\label{luna}
\begin{enumerate}

\item If $W$ is a $G$-module with low separability index, then the
action of $G$ on $W$ is {\it strongly separable} i.e., the stabilizer at
any point is absolutely reduced.  

\item (A version of Luna's \'etale slice theorem in
char.$p$) Let $W$ be a $G$-module with low separability index. Let $F$ be a
fibre of the good quotient $q :W \longrightarrow W//G$, and let $F^{cl}$ be the
unique closed orbit contained in $F$. Then there exists a $G$-map
\[
F \longrightarrow F^{cl}.
\]

\item (see \cite[Proposition 8.5, p.312]{bardsley}) More generally, if $F$ is an affine $G$--subvariety of ${\mathbb P}(W)$, with $W$ as a $G$--module with low separability index, and suppose that $F$ contains a unique closed orbit $F^{cl}$.  Then there exists a $G$-retract
\[
F \longrightarrow F^{cl}.
\]
\end{enumerate}
\eprop

\brem In \cite{bardsley}, Bardsley and Richardson make the assumption that the action of $G$ on $F$ is {\em separable} and the stabilizer $G_f$ at $f \in F^{cl}$ is {\em linearly reductive}. Since the action on $W$ has low separability property, the stabilizer $G_f$ is a saturated, reduced and reductive subgroup of $G$. The assumption of linear reductivity of the stabilizer is handled in our situation by the low separability assumption on $W$, since the tangent space $T_{_f}(F^{cl})$ at $f \in F^{cl}$ is also a $G_f$--module of low separability index. This  gives complete reducibility of the action of the stabilizer on the tangent spaces in Luna's slice theorem (for details see \cite[Proposition 5, p 16]{bapa}). \erem

\bdefe\label{socledef} Let $(V, \theta)$ be a semistable Hitchin pair with $deg(V) = 0$.  The Higgs socle subpair $(\psi(V), \theta_{_{\psi(V)}})$ is defined as the sum of all stable subpairs $(W, \theta_W) \subset (V, \theta)$ of degree $0$.\edefe 

\brem The socle subpair $(\psi(V), \theta_{_{\psi(V)}})$ can be easily seen to be a {\em direct sum} of certain stable subpairs of $(V, \theta)$ each of degree $0$. Moreover, if $(V, \theta)$ is {\em not polystable i.e, it is not direct sum of stable Hitchin pairs}, then $(\psi(V), \theta_{_{\psi(V)}})$ is a {\em proper subpair} of $(V, \theta)$ and conversely. \erem

Recall the notion of admissible reductions and polystability of Hitchin pairs (see Definition \ref{stabilitynotions} (\ref{admis}) and  (\ref{poly})).
 
\brem Let $(V, \theta)$ be a semistable Hitchin pair of rank $n$ with $det(V) \simeq \co_X$.  Let $(E, \theta)$ be the underlying principal $SL(n)$--Hitchin pair.
The Higgs socle $\psi(V) \subset V$ gives an exact sequence of Hitchin pairs:
\beqa\label{soclered}
0 \to \psi(V) \to V \to V/{\psi(V)} \to 0
\eeqa
i.e an admissible Higgs reduction of structure group of  $(E, \theta)$ to a maximal parabolic subgroup $\mathcal P \subset SL(n)$. We will call this reduction $(E_{_\mathcal P}, \theta_{_\mathcal P})$ the ``socle" reduction. Note that  $(E, \theta)$ is a polystable Hitchin pair if and only if the parabolic subgroup $\mathcal P$ coincides with the structure group $ SL(n)$ of $E$. \erem
 
\blem\label{socle} Let $(V, \theta)$ be a semistable Hitchin pair of rank $n$ with $deg(V) = 0$.  Let $(E, \theta)$ be the underlying principal $GL(n)$--Hitchin pair. Let $\mathcal P \subset GL(n)$ be the maximal parabolic subgroup coming from the Higgs socle of $(V, \theta)$ and let $(E_{_\mathcal P}, \theta_{_\mathcal P})$ be the socle Higgs reduction. Then $(V, \theta)$ (or equivalently $(E, \theta)$) is {\em not} polystable if and only if there is no Levi reduction of structure group of $(E_{_\mathcal P}, \theta_{_\mathcal P})$. \elem

{\it Proof}: The proof is trivial since the Higgs socle $(\psi(V), \theta_{_{\psi(V)}})$ is the maximal subpair which has the defining properties in Definition \ref{socledef}. A reduction of structure group of $(E_{_\mathcal P}, \theta_{_\mathcal P})$ to its Levi would mean a splitting of  \eqref{soclered} and this would contradict the maximality of the Higgs socle. \begin{flushright} {\it QED} \end{flushright}

\blem\label{stab1}(cf. \cite[Lemma 10, page 20]{bapa}) Let $(E,\theta)$ be a stable principal $G$--Hitchin pair  with $G$ {\em semisimple}. Suppose that $M$ is a finite dimensional $G$--module with low separability index. Let $Z = SL(M)/Q$ where $Q \subset SL(M)$ is a maximal parabolic subgroup and let $\cl = \cl_{\eta}$ be a very ample line bundle on $Z$ coming from a dominant character $\eta$ of $Q$.  Let $(E(Z), \theta_Z)$ be the associated Hitchin pair. Then any non-zero Higgs section $\sigma: X \lr E(Z)$ such that $deg(s^*{\cl}) = 0$ is a GIT semistable Higgs section in the sense of Definition \ref{gitunstable}.  \elem
 
{\it Proof}: Let $Q \subset SL(W)$ be a  maximal parabolic  and $\eta$  a dominant character of $Q$. Suppose that we are given a Higgs-reduction of structure group $s:X \to E(SL(W)/Q)$. 
{\em Suppose further that it is GIT unstable} in the sense of Definition \ref{gitunstable}. We will get a contradiction to the stability of $(E,\theta)$.
The proof follows almost verbatim the proof of the Theorem \ref{maintheorem} till we reach the construction of the line bundle $L_0$. We pick the thread there. 

Since $(E,\theta)$ is stable, by the Theorem \ref{maintheorem}, the associated vector bundle $E(V_j)$ is Higgs semistable of degree $0$ implying that $deg(L_0) \leq 0$. On the other hand, by \eqref{mtomnought}, we get a map $s^*{\cl} \to L_0$. Since $deg(s^*{\cl}) = 0$, this implies that $deg(L_0) = 0$. 

Now by Proposition \ref{rrkey},  $m_0$ is
in fact GIT semistable for the action of $P_{_\lambda}/U$, i.e for the Levi quotient of $P_{_\lambda}$, for a
suitable choice of linearisation obtained by twisting the action by a
dominant character $\chi$ of $P_{_\lambda}$ (see \eqref{chiandlambda}).  Let $E_{_{P_{_\lambda}}}$ be the $P_{_\lambda}$--bundle obtained by the Higgs reduction to $P_{_\lambda}$.  Let $L_{\chi} = E_{_{P_{_\lambda}}}(\chi)$. 

The GIT semistability of the point $m_0$ for the action of $P_{_\lambda}/U$ with respect to this new linearisation, forces, by Lemma \ref{degandss}, the degree inequality, $deg(L_0 \otimes L_{\chi}^{\vee}) \leq 0$ (see proof of Theorem~\ref{maintheorem}). Since $deg(L_0) = 0$, this in conjunction  with the degree inequality above, gives 
$deg(L_{\chi}^{\vee}) \leq 0$. This implies that for the dual line bundle we have  $deg(L_{\chi}) > 0$. This inequality contradicts the Higgs stability of $(E,\theta)$ (see Definition \ref{stabilitynotions} (i)). \begin{flushright} {\it QED} \end{flushright}

With these  results we now have the following basic theorem on polystability (cf. \cite[Lemma 8.3]{bogo} and \cite[Theorem 3.18]{rr}). Note that we need extra assumptions on characteristic of the ground field for polystability of associated constructions to hold. Let $M$ be a finite dimensional $G$--module. 
Denote by 
\beqa\label{psibar}
{\overline {\psi_{_G}(M)}} = max_{_i}\{\psi_{_G}(\wedge^i(M))\} 
\eeqa

\bth\label{polystablehiggs} Let $(E,\theta)$ be a stable Hitchin pair of {\em degree zero} with $G$ semisimple and $\rho: G \rightarrow SL(M)$,  be a representation such  that $p > {\overline {\psi_{_G}(M)}}$. Then the associated Hitchin pair $(E(M), \theta_M)$ is polystable. \eeth

\noindent {\it Proof}:  By Theorem \ref{maintheorem}, since $p > ht_G(M)$, the extended Hitchin pair $(E(M), \theta_M)$ is {\em semistable}. Let $Q \subsetneq SL(M)$ be a proper maximal parabolic subgroup such that $(E(SL(M)), \theta_M)$ has an admissible Higgs reduction to a parabolic subgroup. In other words, in the language of Higgs bundles, we have a degree zero Higgs subbundle of the Higgs bundle $(E(M), \theta_M)$.

Let us denote the projective variety $SL(M)/Q$ by $Z$ and let $\cl$ be a very ample line bundle with a $SL(M)$--linearization (and hence a $G$--linearization) on $Z$. We use the same notation for the induced line bundle on $E(Z)$ as well. Let $s:X \to  E(Z)$ be  a Higgs section which gives a  Higgs line subbundle $s^*(\cl)$ with $deg(s^*(\cl)) = 0$, i.e the section $``s"$ gives an admissible reduction of structure group.

 By Lemma \ref{stab1}, the section $s:X \to E(Z)$ is {\em a GIT semistable section} in the sense of Definition \ref{gitunstable}. By the last part of Lemma \ref{degandss}, we see that the section $s$ takes its values in $F$, where $F \to Z^{ss} \to Z//G$ {\em is a single GIT fibre}; i.e $s:X \to E(F) \subset E(Z)$. Moreover, $F$ is a $G$--invariant affine variety.

By the assumption of low separability index of the module $W = H^0(SL(M)/Q, \cl)$, it follows that the $G$--action on $F$ is {\em separable}. Furthermore, being a GIT fibre, it contains a unique closed orbit $F^{cl}$. Let $F^{cl} = G/I$. By the low separability assumptions, it follows that the stabilizer $I$  is {\em reduced} and the affineness of $F$ implies that $I$ is reductive.  Moreover, $I \subset Q$.

Again, by virtue of the low separability of the $G$--module $W$, we can apply Proposition \ref{luna}  and we have a $G$--retract, $F \to F^{cl}$, which gives by composition, a section $s':X \to E(F^{cl})$. The section $s':X \to E(F^{cl}) = E(G/I)$, being a Higgs section, gives a Higgs reduction of structure group of the $G$--Hitchin pair $(E,\theta)$ to $I$. Denote this $I$--Hitchin pair by $(E_I, \theta_I)$.

The stabilizer is a {\em saturated subgroup} and hence by the low height property of the representation $\rho:G \to SL(M)$, it follows that the inclusion $I \hra Q$ is {\em completely reducible}. That is, $I \hra L$ for a Levi subgroup of $Q$. (See \cite[Page 15-16]{bapa} and \cite[Page 25]{mour} for details on ``saturated subgroups".)

The upshot of the discussion is that the Hitchin pair $(E(SL(M)), \theta_M)$ gets a Higgs reduction of structure group to the Levi subgroup $L \subset Q$. The argument has shown that whenever $Q \subsetneq SL(M)$ gives an admissible reduction of structure group of 
$(E(SL(M)), \theta_M)$, it gets a further Levi reduction. Coupled with Lemma \ref{socle}, we see immediately that, since the socle reduction $(E_{_\mathcal P}, \theta_{_\mathcal P})$ is an admissible reduction, the parabolic subgroup $\mathcal P \subset SL(M)$ corresponding to the Higgs ``socle"   cannot possibly be a proper parabolic, i.e $\mathcal P = SL(M)$. This implies that $(E(M), \theta_M)$ is a polystable $SL(M)$--Hitchin pair. This completes the proof of the theorem.
\begin{flushright} {\it QED} \end{flushright}

\brem The above theorem gives by far the best effective bounds on $p$ for polystability of associated bundles to hold, even in the situation when there are no Higgs structures on the bundles. \erem

\brem\label{correction} In \cite[Theorem 3.1]{paramvikram} (cf. also \cite[Section 4]{mehtaicm}) it was claimed that if $p$ is larger than the maximum of the {\em heights} of exterior powers as in \eqref{psibar} (and not the separability index) then it gives the polystability of the associated bundle. The proof of \cite[Theorem 3.1]{paramvikram}, which was needed to justify this claim, is incorrect as was pointed out to the authors of \cite{paramvikram} by Professor J. P. Serre in a private correspondence. \erem

\brem The previous theorem is the precise algebraic counterpart of the differential geometric fact that a polystable Higgs vector bundle supports a uniquely defined Einstein-Hermitian connection, a fact proven by Hitchin over curves and by Simpson for smooth projective varieties. \erem

\brem We believe that the approach in this paper should generalize to other natural situations such as ``quiver bundles". Suitable analogues of Hitchin schemes defined for them along with the GIT developed in this paper should yield similar theorems. \erem

\brem\label{higherdim} ({\it Rational Principal bundles}) Following (\cite[page 290]{rr}), we have the corresponding notion of a rational principal $G$--Hitchin pair {\em on a higher dimensional smooth projective variety} $X$. It goes without saying that since we deal with $\mu$--semistability and stability, all the results proven in the previous sections go through without change for $\mu$--semistable (resp. $\mu$--stable) rational $G$--Hitchin pairs. \erem

\end{document}